\documentclass[reqno,12pt]{amsart}

\usepackage{graphics, color, amsmath, amsfonts, amssymb, amsthm, amscd}
\usepackage{graphicx}

\usepackage{epsf}

\def\writefig#1 #2 #3 {\rlap{\kern 62.5 truemm\kern #1 truecm

    \raise #2 truecm \hbox{#3}}}

\usepackage{psfrag}

\usepackage{hyperref}

\numberwithin{equation}{section}

\newcommand{\red}{\textcolor{black}}
\newcommand{\green}{\textcolor{black}}

\DeclareMathSymbol{\leqslant}{\mathalpha}{AMSa}{"36} 
\DeclareMathSymbol{\geqslant}{\mathalpha}{AMSa}{"3E} 
\DeclareMathSymbol{\eset}{\mathalpha}{AMSb}{"3F}     
\renewcommand{\leq}{\;\leqslant\;}                   
\renewcommand{\geq}{\;\geqslant\;}                   
\DeclareMathOperator*{\union}{\bigcup}       

\newcommand{\sumtwo}[2]{\sum_{\substack{#1 \\ #2}}} 
\newcommand{\sumthree}[3]{\sum_{\substack{#1 \\ #2 \\ #3}}} 
\newcommand{\uniontwo}[2]{\union_{\substack{#1 \\ #2}}} 

\makeatletter
\def\captionfont@{\footnotesize}
\def\captionheadfont@{\scshape}
\long\def\@makecaption#1#2{%
  \vspace{2mm}
  \setbox\@tempboxa\vbox{\color@setgroup
    \advance\hsize-6pc\noindent
    \captionfont@\captionheadfont@#1\@xp\@ifnotempty\@xp
        {\@cdr#2\@nil}{.\captionfont@\upshape\enspace#2}%
    \unskip\kern-6pc\par
    \global\setbox\@ne\lastbox\color@endgroup}%
  \ifhbox\@ne 
    \setbox\@ne\hbox{\unhbox\@ne\unskip\unskip\unpenalty\unkern}%
  \fi
  \ifdim\wd\@tempboxa=\z@ 
    \setbox\@ne\hbox to\columnwidth{\hss\kern-6pc\box\@ne\hss}%
  \else 
    \setbox\@ne\vbox{\unvbox\@tempboxa\parskip\z@skip
        \noindent\unhbox\@ne\advance\hsize-6pc\par}%
\fi
  \ifnum\@tempcnta<64 
    \addvspace\abovecaptionskip
    \moveright 3pc\box\@ne
  \else 
    \moveright 3pc\box\@ne
\nobreak
\vskip\belowcaptionskip
\fi
\relax
}
\makeatother

\def\writefig#1 #2 #3 {\rlap{\kern #1 truecm
\raise #2 truecm \hbox{#3}}}

\newtheorem{thm}{Theorem}[section]
\newtheorem{lem}[thm]{Lemma}
\newtheorem{prop}[thm]{Proposition}

\newtheorem{rem}[thm]{Remark}

\newcommand{\cA}{\ensuremath{\mathcal A}}

\newcommand{\cC}{\ensuremath{\mathcal C}}

\newcommand{\cE}{\ensuremath{\mathcal E}}
\newcommand{\cF}{\ensuremath{\mathcal F}}

\newcommand{\cH}{\ensuremath{\mathcal H}}
\newcommand{\cI}{\ensuremath{\mathcal I}}

\newcommand{\cL}{\ensuremath{\mathcal L}}

\newcommand{\cN}{\ensuremath{\mathcal N}}

\newcommand{\cR}{\ensuremath{\mathcal R}}

\newcommand{\cU}{\ensuremath{\mathcal U}}
\newcommand{\cV}{\ensuremath{\mathcal V}}

\newcommand{\frC}{\ensuremath{\mathfrak C}}

\newcommand{\frE}{\ensuremath{\mathfrak E}}

\newcommand{\frW}{\ensuremath{\mathfrak W}}

\newcommand{\frp}{\ensuremath{\mathfrak p}}

\newcommand{\bbB}{{\ensuremath{\mathbb B}} }

\newcommand{\bbE}{{\ensuremath{\mathbb E}} }

\newcommand{\bbN}{{\ensuremath{\mathbb N}} }

\newcommand{\bbP}{{\ensuremath{\mathbb P}} }
\newcommand{\bbQ}{{\ensuremath{\mathbb Q}} }
\newcommand{\bbR}{{\ensuremath{\mathbb R}} }
\newcommand{\bbS}{{\ensuremath{\mathbb S}} }

\newcommand{\bbZ}{{\ensuremath{\mathbb Z}} }

\newcommand{\gep}{\varepsilon}       

\newcommand{\gl}{\lambda}

\newcommand{\gs}{\sigma}

\newcommand{\sfn}{{\sf n}}
\newcommand{\sfx}{{\sf x}}
\newcommand{\sfy}{{\sf y}}
\newcommand{\sfz}{{\sf z}}
\newcommand{\sfu}{{\sf u}}
\newcommand{\sfv}{{\sf v}}
\newcommand{\sfw}{{\sf w}}
\newcommand{\sft}{{\sf t}}

\newcommand{\sfe}{{\sf e}}

\newcommand{\Cl}{\mathbf{Cl}}

\newcommand{\Kp}{{\bf K}_p}
\newcommand{\pKp}{\partial\Kp}

\newcommand{\Perc}{\bbB_p}           

\def\1{\ifmmode {1\hskip -3pt \rm{I}}
\else {\hbox {$1\hskip -3pt \rm{I}$}}\fi} 

\newcommand{\smallo}{o}
\newcommand{\so}{\smallo (1)}
\newcommand{\df}{\stackrel{\Delta}{=}}
\newcommand{\nnsim}[1]{\stackrel{#1}{\sim}}
\newcommand{\eqvs}{\stackrel{\sim}{=}}
\newcommand{\sep}{~\Big|~}

\newcommand{\abs}[1]{\lvert#1\rvert}  

\newcommand{\setof}[2]{\left\{#1\,:\,#2\right\}}
\newcommand{\bigsetof}[2]{\Bigl\{#1\,:\,#2\Bigr\}}

\newcommand{\lb}{\left(}
\newcommand{\rb}{\right)}
\newcommand{\lbr}{\left\{}
\newcommand{\rbr}{\right\}}
\newcommand{\la}{\left\langle}
\newcommand{\rab}{\right\rangle}
\newcommand{\ran}[2]{\right\rangle^{#1}_{#2}} 

\newcommand{\lra}{\longleftrightarrow}
\newcommand{\lraf}{\stackrel{{\rm f}}{\lra}}

\newcommand{\slra}[1]{\stackrel{{#1}}{\lra}}

\newtheorem{bigthm}{Theorem}   

\newcommand{\wtilde}{\widetilde}     

\newcommand{\leqs}{\lesssim}            
\newcommand{\geqs}{\gtrsim}             

\newcommand{\wt}{\widetilde}

\newcommand{\sas}[1]{{#1}^{+}}   
 \newcommand{\sao}[1]{{#1}^{0}}                      

\newcommand{\lsp}{\left [}           
\newcommand{\rsp}{\right ]}          

\newcommand{\rabs}{\right |}
\newcommand{\lpr}{\left .}

\newcommand{\be}{\begin{equation}}
\newcommand{\ee}{\end{equation}}

\begin{document}

\title[\,]{
Finite Connections
 for Supercritical Bernoulli Bond Percolation in 2D
}






\author{Massimo Campanino}
\address{
Dipartimento di Matematica \\
Universit\`a di Bologna \\
piazza di Porta S. Donato 5, I-40126 \\
Bologna, Italy}
\email{campanin\@@dm.unibo.it}
\thanks{}

\author{Dmitry Ioffe}
\address{
Faculty of Industrial Engineering\\
Technion, Haifa 3200, Israel}
\email{ieioffe\@@ie.technion.ac.il}
\thanks{Partly supported by the Japanese Technion Society}

\author{Oren Louidor}
\address{
  Courant Institute of the Mathematical Sciences\\
  New York University\\
  251 Mercer Street\\
  New York, NY 10012, USA}
\email{oren.louidor\@@gmail.com}
\thanks{The research of O. Louidor was supported in part by U.S. NSF
grants DMS-0606696 and OISE-0730136}
\vskip 0.2in


\setcounter{page}{1}






\begin{abstract}
Two vertices $\sfx$ and $\sfy$ are said to be finitely connected if they belong to the
same cluster and this cluster is finite.
We derive sharp asymptotics
\eqref{MasterF2}
of finite connections for super-critical Bernoulli bond
percolation on $\bbZ^2$.
\end{abstract}

\maketitle

\tableofcontents

\section{Introduction and Results}
In the case of the two dimensional nearest neighbour
Ising model below critical temperature, truncated two-point functions could be
computed exactly,
\begin{equation}
\label{MasterF1}
g (\sfx) = \la \sigma_0 ;\sigma_{\sfx}\ran{\,}{\beta} \, =\,
\frac{\phi (\sfn_{\sfx})}{\abs{\sfx}^2}
{\rm e}^{-2\tau_\beta (\sfx)}
\lb 1 +\so\rb ,
\end{equation}
where $\tau_\beta$ is the surface tension,  $\sfn_{\sfx}= \sfx/\abs{\sfx}\in \bbS^1$ and
 $\phi$ is a positive locally analytic function on $\bbS^1$.
\medskip

In this paper we rigorously derive a version of
\eqref{MasterF1} for the simplest  non exactly solvable two
dimensional model: the super-critical Bernoulli bond percolation
on two-dimensional square lattice. The model is self-dual: let $p^* >1/2$,
and consider
 sub-critical  Bernoulli bond percolation measure $\Perc$ on
the direct lattice
$\bbZ^2$ with $p=1-p^*$. Let $\cE^2$ be the set of all nearest neighbour
direct bonds.  Each
direct bond $b\in \cE^2$ intersects exactly one
dual bond $b_*\in \cE^2_*$ of the dual lattice $\bbZ_*^2 = (1/2 ,1/2)+\bbZ^2$.
 Thus each direct percolation configuration
 $\eta\in \lbr 0,1\rbr^{\cE^2}$ unambiguously corresponds to the
dual configuration $\eta_* \in \lbr 0,1\rbr^{\cE^2_*}$ via
\[
\eta (b)\, =\, 1\quad \Longleftrightarrow\quad \eta_* (b_*)\, =\, 0 .
\]
Of course, the induced measure on $\lbr 0,1\rbr^{\cE^2_*}$ is
just the super-critical Bernoulli bond percolation at $p^*$ and we shall
causally take advantage of the fact that both models are
defined on the same probability space and, furthermore, we shall use the
same notation $\Perc$ for both.

Two dual lattice
points $\sfx^* ,\sfy^*\in\bbZ^2$ are said to be finitely connected;
$\{ \sfx^*\lraf \sfy^*\}$, if there
exists  a path of open dual bonds $\gamma^*$ leading
from $\sfx^*$ to $\sfy^*$, but the cluster $ \Cl (\sfx^* )$ of
$\sfx^*$ (and hence
the cluster $\Cl (\sfy^* )$) is finite. The truncated
two-point function is defined then as
\[
 g (\sfx^* ) = \Perc \lb 0^*\lraf \sfx^*\rb ,
\]
where  $0^* \df (1/2, 1/2)$. For simplicity we
shall consider only on-axis directions, that
is we shall focus on asymptotics of $g (\sfx_N^* )$ for
$\sfx_N^* \df  (N+1/2 ,1/2 )$. It should be noted, however, that our
approach goes through with only minor modifications for arbitrary
lattice directions.

\begin{bigthm}
\label{bigthm:A}
For every $p^* =1-p >p_c =1/2$ there exists a constant $\psi =\psi (p^*) >0$
such that
\begin{equation}
\label{MasterF2}
g (\sfx_N^* )\, =\, \Perc\lb 0^*\lraf \sfx_N^*\rb\, = \, \frac{\psi}{N^2}
{\rm e}^{-2N\tau_p (\sfe_1 )} \lb 1+ \so\rb ,
\end{equation}
where $\sfe_1 = (1,0)$ and
 $\tau_p(\cdot)$ is the inverse correlation length of the sub-critical
model (equivalently, the surface tension of the dual  super-critical
model).
\end{bigthm}
The logarithmic asymptotics ${\rm e}^{-2N\tau_p (\sfe_1 )}$ can be
established by relatively soft arguments \cite{CCGKS}: roughly speaking
the event
$ \lbr 0^*\lraf \sfx_N^*\rbr$ implies two disjoint sub-critical connections
over the strip $\lbr (x,y) : 0\leq x\leq N\rbr$. The main struggle here is to recover
asymptotics of finite connection probabilities up to zero order terms, the correct
order $1/N^2$ of the prefactor in particular. This amounts to developing a detailed stochastic geometric characterization of long finite super-critical clusters, which
may be considered as the principle new result of this paper.

Sharp asymptotitcs of finite connections for $d$-dimensional ($d\geq 3$) high-density
models were recently investigated in \cite{BPS1,BPS2}. In the case of higher dimensions
the order of the prefactor is $N^{-(d-1)/2}$.
This is the classical off-critical Ornstein-Zernike prefactor.
The expected non-Ornstein-Zernike order
of prefactor $N^{-2}$ in \eqref{MasterF2} in two-dimensions
 was clearly
understood and discussed on heuristic level
in an earlier literature, see e.g. \cite{BF,BPS2}: The conventional OZ
picture comes from the fluctuation theory of one-dimensional systems. However, finite connections
in two dimensions are described in terms of fluctuation theory of
 two interacting one dimensional
effective random walk type structures.

Let us elaborate on the latter point.
Both the direct
sub-critical percolation model at $p<1/2$ and the dual super-critical
model at $p^* = 1-p>1/2$  are
defined on the same probability space.

In particular, the event $\lbr 0^*\lraf \sfx_N^*\rbr$ can be
written  {as} (see Figure~1)
\begin{equation}
\label{Dualform}
\lbr 0^*\lraf \sfx_N^*\rbr\, =\, \lbr 0^*\lra \sfx_N^*\rbr\cap \cC_N ,
\end{equation}
where $\lbr 0^*\lra \sfx_N^*\rbr$ means that $0^*$ and $\sfx_N^*$ are connected
in the dual  model and the event $\cC_N$ is defined in terms of
the direct  percolation model via
\begin{equation}
\label{EventCN}
\cC_N \, =\ \, \bigsetof{\eta \in \lbr 0,1\rbr^{\cE^2}}{\exists\,\,
\text{an open direct  loop around $0^*$ and $\sfx_N^*$}} .
\end{equation}

\begin{figure}[tbh] 
\psfrag{a}{$0^*$}
\psfrag{b}{$\sfx_N^*$}
\psfrag{c}{$\gamma^*$}
\psfrag{d}{$\frC_N$}
\psfrag{x}{$\sfx$}
\psfrag{y}{$\sfy$}
\psfrag{v}{$\sfv$}
\psfrag{u}{$\sfu$}
\psfrag{l}{$\ell$}
\psfrag{N}{$N-r$}
\centerline{\includegraphics{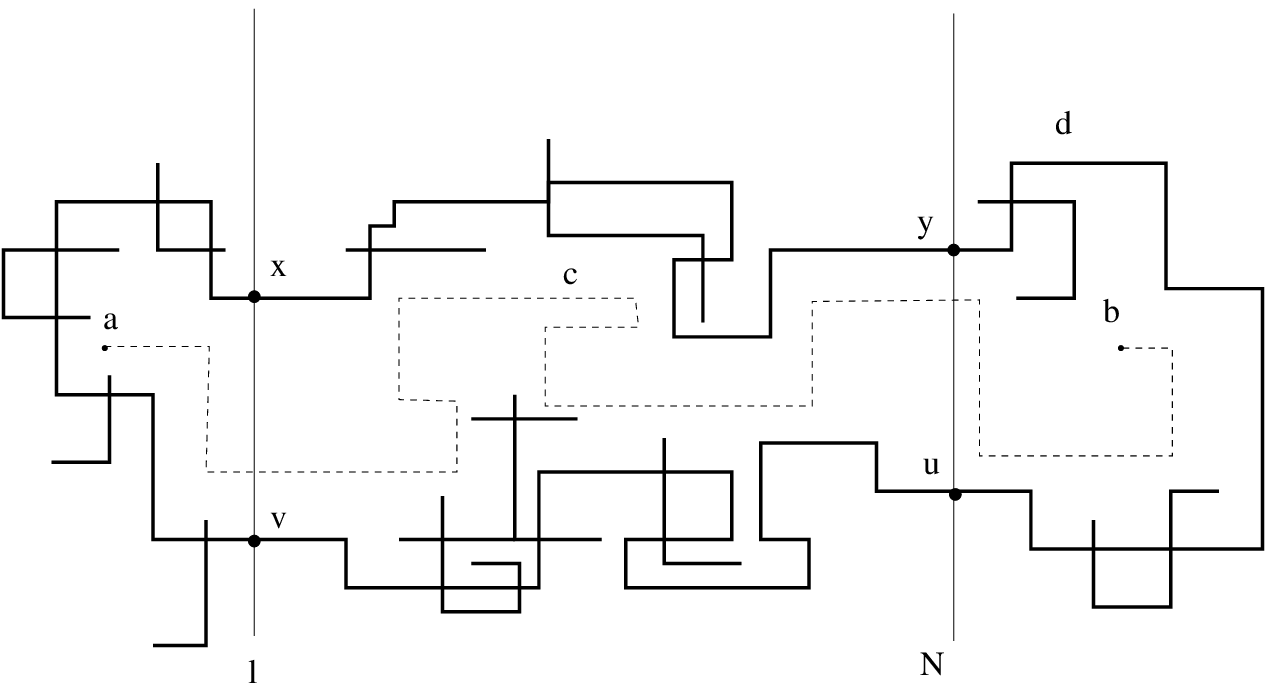}}
\epsfxsize=125mm     
\label{Fig1}
\caption{The event $\lbr 0^*\lraf \sfx_N^*\rbr\,
=\, \lbr 0^*\lra \sfx_N^*\rbr\cap \cC_N$: $\gamma^* :0^*\rightarrow \sfx_N^*$
is
a dual (super-critical) open path, whereas $\frC_N$ is an open direct   (sub-critical) loop-like cluster. The cluster $\frC_N$ splits into
irreducible loops around $0^*$ and $\sfx_N^*$ and a pair of disjoint connections
from $\sfx$ to $\sfy$ and from $\sfv$ to $\sfu$.}

\end{figure}

We shall use $\frC_N$ to denote the inner-most connected component
which contains such a loop, and
we shall decompose $\{0^*\lraf \sfx_N^*\}$ according to
geometric properties of $\frC_N$. We shall see that a typical $\frC_N$
can be split as it is schematically depicted on Figure~1: There are
``irreducible'' dual percolation loops from $\sfv$ to $\sfx$ around $0^*$ and,
respectively, from $\sfu$ to $\sfy$ around $\sfx_N^*$.  In the middle strip
$\lbr (x, y): \ell\leq x\leq N-r\rbr$ there are disjoint connections from $\sfv$ to
$\sfu$ and from $\sfx$ to $\sfy$.  The notion of irreducibility will be set up in
such a way that $\ell$ and $r$ will be typically small and will be
eventually integrated out. The crux of the matter is to understand how to
compute the probability of the double connection event in the middle strip.
The main thrust of the theory developed in \cite{CI,CIV1,CIV3} is that on large finite
scales sub-critical percolation clusters have effective random walk
structure. One of our two main results here is a reformulation of the double
connection event in the middle strip in terms of hitting probabilities
for two effective random walks conditioned on non-intersection. The second
main result is an adjustment of the fluctuation identities introduced in
\cite{AD,BJD} for computing these probabilities.

\subsubsection{Effective random walk picture}

\label{sub:effective_rw}

We  proceed with a description of our effective random walk
picture as it will show up in the reformulation of the double connection
event.
Let $\lbr \sigma_k = (\rho_k, \xi^1_k ,\xi^2_k )\rbr$ be
a collection of
i.i.d. $\bbN\times \bbZ^2$-valued random variables defined on some probability
space equipped with a probability measure $\bbP$
and satisfy the following
set of conditions:

\smallskip

\noindent
{\bf (P1)} There exists $\alpha <\infty $ such that
\[
{\rm Range}(\sigma )\, =\,
\setof{(t, v ,x )}{\abs{v},\abs{x}  <\alpha t} .
\]

\noindent
{\bf (P2)} There exists $\beta >0$ such that (for  the exact definition of ``$\leqs$''  see the remark on notational convention below)
\[
\bbP \lb \rho  > t\rb\, \leqs \, {\rm e}^{-\beta t} .
\]

\noindent
{\bf (P3)} For any $t\in\bbN$ the conditional (on $\rho  = t$) distribution
of $(\xi^1 ,\xi^2 )$ is symmetric in $\bbZ^2$ with respect to the axes and the diagonal
$\{(v,x) :\; v = x\}$, that is for any $(v, x)\in\bbZ^2$:
\begin{equation}
 \begin{split}
\bbP\lb \xi^1 =v, \xi^2 = x\sep \rho =t\rb \, &=\,
\bbP\lb \xi^1 =\abs{v}, \xi^2 = \abs{x}\sep \rho =t\rb\,\\
&  =\,
\bbP\lb \xi^1 =x, \xi^2 = v\sep \rho =t\rb .
\end{split}
\end{equation}
\smallskip

Define random walk $S_n = \lb T_n , V_n ,X_n  \rb \, =\, \sum_1^n\sigma_k +S_0$,
and let $\bbP_{v , x}$ be
the law of this random walk subject to the initial condition
$S_0 = (0,v,x)$.
Consider the following event
\begin{equation}
\label{eq:Rplus}
\cR_n^+\, =\, \lbr X_k > V_k \quad\text{for}\ k=1, \dots, n\rbr .
\end{equation}

\begin{bigthm}
\label{bigthm:B}
There exists a function $U :\bbN\to\bbR_+$
of an at most linear growth, $U (z)\leqs z$,
such that
\begin{equation}
\label{eq:RWRplus}
\bbP_{v , x}\lb
\bigcup_n \lb\lbr S_n = (N,  u , y ) \rbr \cap \cR_n^+\rb\rb\,
\sim\, \frac{U(x -v )U(y  - u )}{N^2} ,
\end{equation}
uniformly in  $v < x$ and $u <y$
 satisfying
$\abs{v}, \abs{x}, \abs{u}, \abs{y} \leqs \log N$.
\end{bigthm}

The above function $U$ is in fact a  certain
renewal  function  related to the differences process $Z_n = X_n - V_n$, which is again
a random walk.

\noindent
Theorem~\ref{bigthm:B} follows by an adjustment of one-dimensional techniques
developed in \cite{AD,BJD}. It should be noted that \eqref{eq:RWRplus} could be
extended to a much larger range of parameters $N,v,x,u$ and $y$.

\subsubsection{Organization of the paper.}
The paper is organized as follows: In Section~\ref{sec:TopLevelDecom} we describe
the percolation geometry of finite connections.
We start by introducing the $\bbZ^2$-lattice notation and by recalling
the results of \cite{CI, CIV3} on the geometry of long sub-critical clusters.
 This sets up the stage for basic geometric decomposition
%
\eqref{eq:decomposition} and
\eqref{eq:decompositionT}  of
$\Perc\lb 0^* \lraf \sfx_N^*\rb$. Main claims behind the proof of
Theorem~\ref{bigthm:A} are collected in Subsection~\ref{sub:ProofA}.
As it is explained in Subsections~\ref{sub:structure1}-\ref{sub:Proofdecomposition} both the
the validity of \eqref{eq:decomposition} and the claim of Lemma~\ref{lem:LR} follow by
a more or less straightforward adjustment of the techniques developed in \cite{CI,CIV3}.

The crucial point is to prove Theorem~\ref{lem:A}. The proof is based on the
effective random walk representation \eqref{eq:ARplus} and it is explained in
Subsection~\ref{sub:proofA}. Apart from justifying and establishing various properties
of the representation in Section~\ref{sec:Reduction} there are two types of results involved:
We need a certain generalization (Theorem~\ref{bigthm:C})  of the results of \cite{AD,BJD}
 on random walks
conditioned to stay positive. This issue is addressed in Section~\ref{sec:rw}.
In the concluding Section~\ref{sec:Decoupling} we
develope  estimates on repulsion of effective random walk trajectories and on
decoupling of the associated percolation events.
\medskip

\subsubsection{Remark on notational conventions}
Let $\{ a_n ( w)\}$ and $\{ b_n (w )\}$ be two sequences of positive
numbers indexed by $w$ from some set of parameters $w\in\frW_n$.
We say that $a_n (w)\sim b_n (w) $
if there exists
a constant $c >0$, such that
\[
\lim_{n\to\infty}\frac{a_n (w)}{b_n (w)}\, =\, c
\]
uniformly in $w\in\frW_n$.
If we want to specify the exact value of the  constant $c$
appearing above, we shall write $a_n (w)\nnsim{c} b_n (w) $

Similarly, let us say that $a_n ( w)\leqs b_n (w )$, if
\[
\limsup_{n\to\infty}\frac{a_n (w)}{b_n (w)}\, <\, \infty ,
\]
uniformly in $w\in\frW_n$. Often the dependence on $w$ will
not be written explictly and furthermore, in some cases, there
will be no additional parameter at all. Where confusion arises, we shall
indicate the dependency (or lack of it) explicitly. In addition, that same notation
will be used to specify $\frW_n$ itself. For example, if we say that a certain
property holds uniformly in $|w| \leqs v_n$, where $v_n$ is a given sequence,
then for every $K$ fixed this property holds if $|w| \leq K v_n$
and $n$ is large enough.
Finally, let us say that $a_n (w )\eqvs b_n (w)$  if there exists constant
$c$ such that $a_n (w ) =c b_n (w)$ for
all $w\in \frW_n$.

In the sequel we shall often rely on the following relation,
 which we call {\em Gaussian summation formula}: Let $A$ be a non-degenerate
quadratic form on $\bbR^d$. Then,
\[
 \sum_{x\in\bbZ^d} {\rm e}^{-A (x)/n}\, \sim \, n^{d/2} .
\]





\section{Geometry of Finite Connections}
\label{sec:TopLevelDecom}
\subsection{Lattice and dual lattice  notation} 
Most of the work will be done on the direct lattice $\bbZ^2$.
We shall use sans-serif font, e.g. $\sfx, \sfy, \sfu, \dots$  for the 
vertices of $\bbZ^2$ and points in $\bbR^2$ and usual roman font 
to denote their one-dimensional coordinates, e.g. $\sfx = (t, x)$. $|\cdot|$ will 
denote both the absolute values for scalars and the Euclidean norm for vectors.

All quantities which live on the dual lattice $\bbZ^2_*$ are 
marked with $*$, e.g. $\sfx^*$ for vertices, $\sfe^*$ for bonds and 
$\gamma^*$ for paths. 
For each point $\sfx\in\bbZ^2$ define its four ``geographic''  
dual neighbours:
\[
\sfx^*_{ne}= \sfx+ (1/2, 1/2), \ \ 
\sfx^*_{se}= \sfx+ (1/2, -1/2), \ \
 \sfx^*_{sw}= \sfx+ (- 1/2, -1/2), \ \
\sfx^*_{nw}= \sfx+ (-1/2, 1/2) .
\]
Also given  a set $B \subseteq \cE^2$, the set $B^*$ contains all the bonds which are 
 dual to the  bonds in $B$.   

Next define:
\begin{equation*}
\begin{split}
&\cH_m^- = \setof{\sfx =(k,l)\in\bbZ^2}{k< m}\ \ 
 \cH_m^+= \setof{\sfx = (k,l)\in\bbZ^2}{k>  m}\\ 
&\qquad \text{and}\ 
\cH_{m,r} = \setof{\sfx = (k,l)\in\bbZ^2}{m\leq k\leq r}.
\end{split}
\end{equation*}
We shall write $\cH_m$ instead of $\cH_{m,m}$.  
The sets of bonds we associate
with $\cH_{m, r}$ are:
\[
\cE_{m, r}\  =\  \setof{(\sfx ,\sfy)\in\cE^2}{\sfx\in\cH_{m,r}\ \textbf{and}\ 
\sfy\in\cH_{m,r}}\ \ 
\]
and
\[
\cE^+_{m, r} = \cE_{m, r} \setminus \cE_{m, m} \ , \quad 
\cE^-_{m, r} = \cE_{m, r} \setminus \cE_{r, r} \ , \quad
\]
As a shorthand, we write $\cE^+_r$ and $\cE^-_m$ for $\cE^+_{r, \infty}$
and $\cE^-_{-\infty, m}$.
Note that for each $m\leq r$, both $\bbZ^2$ and $\cE^2$ could be represented
as disjoint unions, 
\[
\bbZ^2\, =\, \cH_m^-\vee\cH_{m, r}\vee \cH^+_r\ \ \text{and}\ \  
\cE^2\, =\, \cE_m^-\vee\cE_{ m, r}\vee \cE^+_r .
\]
Let $\frE_{m,r}$ be the $\sigma$-algebra generated by the
direct percolation configuration on $\cE_{m, r}$ and define
$\frE^{\pm}_{m,r}$, $\frE^{\pm}_r$ in an analogous way.
Under $\Perc$, $\frE_m^-$, $\frE_{m ,r}$ and $\frE_r^+$ 
are independent.

Given a set $A\subseteq\cE^2$ and a percolation configuration
$\eta\in \{ 0, 1\}^A$, let us say that $\sfx\slra{A}\sfy$ if 
$\sfx$ and $\sfy$ are connected by a path of open bonds in $\eta$.
Given $m <r$ and a site $\sfx\in\cH_{m ,r} $ let us define $\Cl_{m,r} (\sfx )$
to be the cluster of sites which are connected to $\sfx$ by direct open bonds
in $\cE_{m,r}$. This is a sub-graph of $\lb \cH_{m,r}, \cE_{m,r} \rb$ but we shall
frequently treat it as a subset of bonds or vertices only. For example, for $A\subseteq \cH_{m,r}$ or
$B\subseteq \cE_{m,r}$, we may write $\{\Cl_{m,r} (\sfx )=A\}$ or
$\{\Cl_{m,r}(\sfx )=B\}$ to indicate which sites or bonds comprise the cluster.
Note that an event defined with either of the two conditions, belongs to $\frE_{m,r}$. 

We use $\Cl_{m,r}(\sfx ,\sfy )= \Cl_{m,r} (\sfx )\cap \Cl_{m,r} (\sfy ) $ to denote 
the common cluster of $\sfx, \sfy \in \cH_{m,r}$ inside the strip $\cH_{m,r}$ .
Similarly, we use 
$\Cl^{\pm}_{m,r}(\dots)$ and $\Cl^{\pm}_m(\dots)$ for the corresponging clusters 
restricted to open bonds from 
$\cE^{\pm}_{m,r}$ and $\cE^{\pm}_m$.

Finally $\prec$ stands for the standard lexicographical order on $\bbZ^2$.  That is,  
$(x_1, x_2) = \sfx \prec \sfy = (y_1, y_2)$ if and only if $x_1 < y_1$ or 
$x_1 = y_1\,,\,\, x_2 < y_2$. 

\subsection{Decomposition of $\{0^*\lraf \sfx_N^*\}$ and basic percolation events}
\label{sub:basic}
\label{sub:event}  It is time to describe precisely our basic
geometric decomposition of the event $\lbr 0^*\lraf \sfx_N^*\rbr$ 
(in its representation \eqref{Dualform}) as it was schematically 
depicted on Figure~1. 

 Given $0<l \leq N$ let us say that $\cH_l$ is a {\em cut line}
of $\frC_N$ if the number of points $\# \lb \frC_N\cap\cH_l\rb =2$. Define, 
\[
 \cI(\emptyset )\, \df\, \{0^*\lraf \sfx_N^*\}\cap\lbr \frC_N\ \text{
contains  less than two cut lines}\rbr
\]
In all the remaining cases we can talk about different left-most and right-most cut-lines
of $\frC_N$:
Given $0 <m < N-r \leq N$ and two pairs of points $\sfv , \sfx\in\cH_m$ and
$\sfu ,\sfy\in\cH_{N-r}$ let us say that $\cI ([\sfv , \sfx ], [\sfu ,\sfy ])$
occurs, if
\begin{equation}
\label{eq:cutplane}
\frC_N\cap\cH_m = \lbr \sfv , \sfx\rbr\quad 
\frC_N\cap\cH_{N-r} = \lbr \sfu , \sfy\rbr , 
\end{equation}
but 
\[
\#\lb \frC_N\cap\cH_l\rb >2\  \forall\, l=1,\dots, m-1\quad 
\text{and}\quad
\#\lb \frC_N\cap\cH_{N-l}\rb >2\  \forall\, l=0,\dots, r-1 .
\]
As a result we represent $\{0^*\lraf \sfx_N^*\}$ as the disjoint union
(below 
$\prec$ stands for the lexicographical order relation), 
\begin{equation}
 \label{eq:disjointI}
\{0^*\lraf \sfx_N^*\}\, =\, \bigcup_{0 <m< N-r \leq N}\uniontwo{\sfv\prec \sfx}{\sfv ,\sfx\in \cH_m}
\uniontwo{\sfu\prec \sfy}{\sfu ,\sfy\in \cH_{N-r}}\cI([\sfv , \sfx ], [\sfu ,\sfy ])\, 
\bigcup\, \cI (\emptyset ) , 
\end{equation}
and, accordingly, 
\[
 \label{eq:sumI}
\Perc\lb 
0^*\lraf \sfx_N^*\rb\, =\, 
\sum_{0 <m< N-r \leq N}\sumtwo{\sfv\prec \sfx}{\sfv ,\sfx\in \cH_m}
\sumtwo{\sfu\prec \sfy}{\sfu ,\sfy\in \cH_{N-r}}
\Perc\lb \cI([\sfv , \sfx ], [\sfu ,\sfy ])\rb \, +\, 
\Perc\lb \cI (\emptyset )\rb .
\]
We shall prove that not only $\Perc \lb \cI (\emptyset ) \rb$ is negligible, but in fact one can restrict
attention to events $ \cI([\sfv , \sfx ], [\sfu ,\sfy ])$ with $\sfv ,\sfx$ being sufficiently
close to $0$ and, respectively, $\sfu ,\sfy$ being sufficiently close to $\sfx_N$. Namely, 

\begin{lem}
\label{lem:decomposition}
\begin{equation}
\label{eq:decomposition}
\begin{split}
&\Perc\lb 0^*\lraf \sfx_N^* \rb (1+ \so )\, 
\\
&\, =\, 
\sum_{0 <m< N-r \leq N}\, \, 
\sumthree{\abs{\sfv - 0}, \abs{\sfx - 0} \,\leqs \,\log N}{\sfv\prec \sfx}{\sfv ,\sfx\in \cH_m} \   \ 
\sumthree{\abs{\sfx_N -\sfu}, \abs{\sfx_N- \sfy} \,\leqs\, \log N}
{\sfu\prec \sfy}{\sfu ,\sfy\in \cH_{N-r}}
\Perc\lb \cI([\sfv , \sfx ], [\sfu ,\sfy ])\rb \\
&\qquad\qquad\qquad \qquad\df\, {\sum}_N 
\Perc\lb \cI([\sfv , \sfx ], [\sfu ,\sfy ]) \rb 
\end{split}
\end{equation}
\end{lem}
\noindent
We shall sketch the proof of this lemma in the end of the Section.

\smallskip 

For technical reasons, which
will
become apparent in Lemma~\ref{lem:toproduct} below, it happens to be convenient to 
work with a slight modification $\widetilde{\cI}([\sfv , \sfx ], [\sfu ,\sfy ])$ of 
${\cI}([\sfv , \sfx ], [\sfu ,\sfy ])$, the precise definition is given in \eqref{eq:Itilde}.
 {Before}, we need to introduce a bit of additional notation: Given $m =1, \dots,  {N}$ and 
$\sfw ,\sfz\in\cH_m$, with $\sfw\prec \sfz$, 
 let us say that 
$\Cl_m^- (\sfw ,\sfz )$
is a loop around $0^*$ {\em rooted} at 
$\lb\sfw ,\sfz\rb$ if, 
\be
\label{eq:root} 
\lbr \sfw \slra{\cE_m^-} {\sfz}\rbr,\quad
\lbr\sfw^*_{nw}\slra{ {(\cE_m^-)^*}} 0^*\slra{ {(\cE_m^-)^*}}\sfz^*_{sw}\rbr\quad
\text{and}\quad \lbr \sfw^*_{sw}\slra{ {(\cE_m^-)^*}}\sfz^*_{nw}\rbr .
\ee
There is a completely symmetric definition of rooted loops around  $\sfx_N^*$.

Let $\Cl_m^- (\sfw ,\sfz )$ be a loop around $0^*$ , rooted at $\sfw ,\sfz\in\cH_m$. 
We shall say that $ {1 <} l < m$ is a {\em modified left cut line} of $\Cl_m^- (\sfw ,\sfz )$ if
there exist $\sfv ,\sfx\in\cH_l$ such that, 
\smallskip

 \noindent
{\bf a)} $\Cl_l^- (\sfv ,\sfx)$
is a loop around $0^*$, rooted at $\lb \sfv,\sfx\rb$. 

\noindent
{\bf b)} $\sfx\slra{\cE_{l ,m}^-}\sfz$\quad\text{and}\ \ 
$\sfx =\max\lbr \Cl_{l,m}^- (\sfx ,\sfz )\cap\cH_l\rbr$.

\noindent
{\bf c)} $\sfv\slra{\cE_{l ,m}^-}\sfw$\quad\text{and}\ \ 
$\sfv =\max\lbr \Cl_{l,m}^- (\sfv,\sfw )\cap\cH_l\rbr$.
\smallskip

There is a completely symmetric definition of {\em  modified right cut lines}. 
A loop 
is said to be {\em irreducible} if it does not have  modified cut 
lines.

Let $l$ be a cut line of $\frC_N$  and denote $\lbr\sfw, \sfz\rbr = \frC_N\cap \cH_l$ (with 
$\sfw\prec\sfz$). \red{If $\Cl_l^- (\sfw ,\sfz )$ is not a rooted loop around $0^*$  
then there must exists another disjoint loop
 $\Cl_l^- (\sfu ,\sfv )$ around $0^*$ for some $\sfw\prec\sfu\prec\sfv\prec\sfz$
with $\Cl_l^- (\sfw ,\sfz )\cap \Cl_l^- (\sfu ,\sfv ) = \emptyset$.} 
Indeed it is only in the latter case when the second of 
\eqref{eq:root} is violated. Thus, conditioning on the realizations of 
$\Cl_l^- (\sfw ,\sfz )$ and \red{using the BK inequality, one deduces},
\begin{eqnarray}
\nonumber
\lefteqn{
\log \Perc \lb 
	\begin{array}{l}
		\Cl_l^- (\sfw ,\sfz )\ \text{is not a rooted}\\
		\quad \text{loop around $0^*$}
	\end{array}
\rabs 
\lpr 
		\begin{array}{l}
			\ \ 0^*\lraf \sfx^*_N\\
			\lbr\sfw, \sfz\rbr = \frC_N\cap \cH_l
		\end{array}
\rb} \\
\label{eq:rootBound}
& \leq 	& \red{\log} \lb \sum_{\sfw\prec\sfu\prec\sfv\prec\sfz}
 			\Perc\lb \Cl_l^- (\sfu ,\sfv )\ \text{is a loop around $0^*$} \rb \rb \leqs -l .
\end{eqnarray}

Let us say that a cut line $l$ with $\lbr\sfw, \sfz\rbr = \frC_N\cap \cH_l$ is {\em strong} if 
 both $\Cl_l^- (\sfw ,\sfz )$ and $\Cl_l^+ (\sfw ,\sfz )$ are rooted loops around $0^*$ and, respectively, 
around $\sfx_N^*$.  Inequality \eqref{eq:rootBound} above controls conditional probabilities that 
$l$ is a strong cut line given that it is a cut line.

Note now that if $1<  k < l <  m$ 
 and $l$ is 
a left modified cut line of 
$\Cl_m^- (\sfw ,\sfz )$ with $\lb\sfv ,\sfx\rb$
beeing the corresponding root, then $k$ is a left modified  cut line of 
$\Cl_m^- (\sfw,\sfz )$ if and only if it is a left modified cut line 
of $\Cl_l^- (\sfv ,\sfx )$.  In paricular, 
once $\frC_N$ contains at least one {\em strong}  
cut line the notions of the left-most left modified 
cut line of $\frC_N$ 
and, accordingly, of the  right-most right modified 
 cut lines  of $\frC_N$, are well defined. 

\subsubsection{Events\  $\widetilde{\cI}([\sfv , \sfx ], [\sfu ,\sfy ])$} 
The events $\widetilde{\cI}([\sfv , \sfx ], [\sfu ,\sfy ])$ are defined for 
$0< l<N-r  {\leq} N$;   $\sfv ,\sfx\in\cH_l$ and $\sfu ,\sfy\in\cH_{N-r}$.
 They are defined in such a way that
  they are disjoint for different choices of $\sfv ,\sfx, \sfu$ and $\sfy$. Moreover, 
once $\lbr 0^*\slra{f}\sfx_N^*\rbr$ occurs and $\frC_N$ has at least
two cut lines, one of $\widetilde{\cI}([\sfv , \sfx ], [\sfu ,\sfy ])$ necessarily 
happens. Loosely speaking, the event $ \widetilde{\cI}([\sfv , \sfx ], [\sfu ,\sfy ])$ 
requires that $l$ and $N-r$ are the left-most (respectively right-most)
modified left (respectively right) cut lines with the corresponding irreducible loops being 
rooted at $(\sfv ,\sfx )$ (respectively $(\sfu, \sfy)$). 
 Formally, $\widetilde{\cI}([\sfv , \sfx ], [\sfu ,\sfy ])$ is represented as an 
intersection of three independent events, 
\be
\label{eq:Itilde}
 \widetilde{\cI}([\sfv , \sfx ], [\sfu ,\sfy ])\, =\, 
\cL([\sfv ,\sfx])\cap\cA ([\sfv ,\sfx], [\sfu ,\sfy ])\cap \cR([\sfu ,\sfy ]) .
\ee
\subsubsection{Events $\cL([\sfv ,\sfx])$ and $\cR([\sfu ,\sfy ])$}
For $\sfv ,\sfx\in\cH_l$, the event $\cL([\sfv ,\sfx])$ is defined as 
\[
 \cL([\sfv ,\sfx]) = \lbr  \Cl_{l}^- (\sfv ,\sfx )\ \text{is an irreducible loop aroud $0^*$}\rbr .
\]
For $\sfu,\sfy\in\cH_{N-r}$, the event $\cR([\sfu ,\sfy])$ is defined as 
\[
 \cR([\sfu ,\sfy]) = \lbr  \Cl_{N-r}^+ (\sfu ,\sfy )\  
\text{is an irreducible loop aroud $\sfx_N^*$}\rbr .
\]
{(See Figure~2{\it (ii), (iii)}).}

\subsubsection{Events  $\cA ([\sfv ,\sfx], [\sfu ,\sfy ])$} For each 
$m<N-r$, each pair of vertices $\sfv\prec \sfx $; $\sfv,\sfx\in\cH_m$  and 
each pair of vertices $\sfu\prec \sfy $; $\sfu,\sfy\in\cH_{N-r}$ the 
event  $\cA ([\sfv ,\sfx], [\sfu ,\sfy ])$ is defined by the following
set of conditions (Figure~2{\it (i)})
\smallskip

\noindent
{\bf a)} $\Cl_{m,N-r} (\sfv ,\sfu)\neq \emptyset$. 

\noindent
{\bf b)} $\Cl_{m,N-r} (\sfx ,\sfy)\neq \emptyset$.

\noindent
{\bf c)}  $\sfv = \max\{\Cl_{m,N-r} (\sfv ,\sfu)\cap \cH_m\}$ and 
 $\sfu = \max\{\Cl_{m,N-r} (\sfv ,\sfu)\cap \cH_{N-r}\}$, where the maximum 
is understood in the lexicographical order, e.g.  $\sfv$ has the maximal 
vertical coordinate among all the vertices in 
$\Cl_{m,N-r} (\sfv ,\sfu)\cap \cH_m$.

\noindent
{\bf d)} $\sfx = \max\{\Cl_{m,N-r} (\sfx ,\sfy)\cap \cH_m\}$ and 
 $\sfy = \max\{\Cl_{m,N-r} (\sfx ,\sfy)\cap \cH_{N-r}\}$. 

\noindent 
{\bf e)} $\Cl_{m,N-r} (\sfv ,\sfu)\cap \gamma^{\rm up}(\Cl_{m,N-r} (\sfx ,\sfy ))
 =\emptyset$, 
where $\gamma^{\rm up}(\Cl_{m,N-r} (\sfx ,\sfy ))$ is the upper envelope of the
 cluster $ \Cl_{m,N-r} (\sfx ,\sfy)$.
\smallskip

Notice that conditions {\bf c)} and {\bf e)} imply that 
\begin{equation}
\label{eq:vtou}
\sfv^*_{nw}\slra{{(\cE_{m,N-r})^*}}\sfu^*_{ne} .
\end{equation}
On the other hand, 
condition {\bf e)} by itelf  may seem redundant: Indeed in view of the strict ordering 
$\sfv\prec \sfx$ conditions {\bf a)}-{\bf d)} already ensures  that 
 $ \Cl_{m,N-r} (\sfx ,\sfy)\cap\Cl_{m,N-r} (\sfv ,\sfu) =\emptyset$. The
reason for choosing such a formulation will become apparent in Lemma~\ref{lem:toproduct}.

\smallskip

We stress  that 
$\widetilde{\cI}([\sfv , \sfx ], [\sfu ,\sfy ])$ are 
disjoint for different choices of $\sfv,\sfx ,\sfu$ and $\sfy$
{and 
\[
 \lbr 0^*\slra{f}\sfx^*_N\rbr
\supseteq\bigcup_{\sfx ,\sfv}\bigcup_{\sfu ,\sfy}
\widetilde{\cI}([\sfv , \sfx ], [\sfu ,\sfy ]) .
\]
}
Since left-most and right-most modified cut lines
are well defined and distinct whenever $\frC_N$ has at least two strong  cut lines, 
it is rather straightforward to deduce from \eqref{eq:decomposition} and 
 \eqref{eq:rootBound}
that, 
\begin{equation}
\label{eq:decompositionT}
\begin{split}
\Perc\lb 0^*\lraf \sfx_N^* \rb (1+ \so )\, 
&
\, =\, 
{\sum}_N 
\Perc\lb \widetilde\cI([\sfv , \sfx ], [\sfu ,\sfy ])\rb \\
&\, =\, 
{\sum}_N
\Perc\lb \cL([\sfv , \sfx])\rb \Perc\lb \cA([\sfv ,\sfx], [\sfu ,\sfy ])\rb 
\Perc \lb \cR([\sfu ,\sfy ])\rb .
\end{split}
\end{equation}

\begin{figure}[tbh] 
\psfrag{a}{$\sfv^*_{nw}$}
\psfrag{b}{$\sfu^*_{ne}$}
\psfrag{c}{$\sfv$}
\psfrag{d}{$\sfu$}
\psfrag{g}{$\sfx$}
\psfrag{e}{$\gamma$}
\psfrag{f}{$\sfy$}
\psfrag{h}{$\sfy$}
\psfrag{i}{$\sfu$}
\psfrag{j}{$\sfu^*_{ne}$}
\psfrag{k}{$\sfx_N^*$}
\psfrag{l}{$0^*$}
\psfrag{m}{$\sfv^*_{nw}$}
\psfrag{n}{$\sfv$}
\psfrag{o}{$\sfx$}
\psfrag{p}{$\cH_m$}
\psfrag{s}{$\cH_{N-r}$}
\psfrag{r}{$\Cl_{m, N-r}(\sfx ,\sfy )$}
\psfrag{t}{$\Cl_{m, N-r}(\sfv ,\sfu )$}
\psfrag{v}{$\Cl_-$}
\psfrag{u}{$\Cl_+ $}
\centerline{\includegraphics{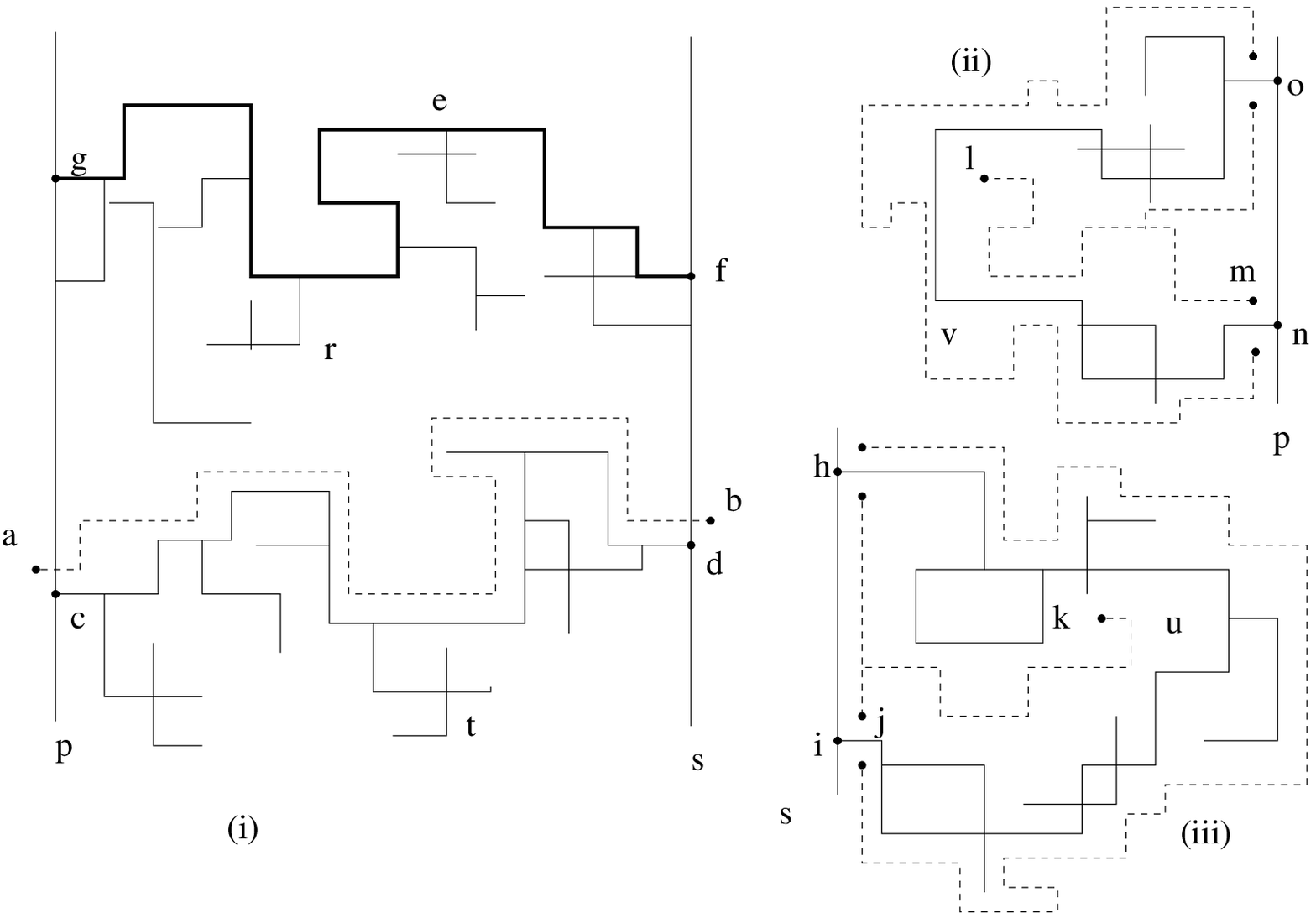}}
\label{Fig2}
\caption{
{\it (i)} Event $\cA([\sfv ,\sfx], [\sfu ,\sfy ])$: $\gamma$ is the upper 
envelope of $\Cl_{m,N-r} (\sfx ,\sfy )$
\newline  {\it (ii)} Event $\cL ([\sfv , \sfx ])$, $\Cl_- \equiv  \Cl_{m}^- (\sfv ,\sfx )$  \quad 
{\it (iii)}  Event $\cR ([\sfu ,\sfv ])$, $\Cl_+\equiv \Cl_{N-r}^+ (\sfu ,\sfy )$ .
}
\end{figure}

\subsection{Proof of Theorem~\ref{bigthm:A}}
\label{sub:ProofA}

The proof of Theorem ~\ref{bigthm:A} will follow immediately from 
\ref{eq:decompositionT}, once we establish Lemma ~\ref{lem:decomposition}, Theorem~\ref{lem:A} 
and 
Lemma~\ref{lem:LR} below.

Recall that $\tau_p(\cdot$) is the inverse correlation length for the 
sub-critical model. Set $\tau_p = \tau_p (\sfe_1 )$ and
$\sft_p = \tau_p \sfe_1  = (\tau_p , 0)$. 
Let us use $\la\cdot , \cdot\rab$ to denote the scalar product in $\bbR^2$. 
Notice that in view of lattice symmetries, 
\[
 \tau_p (\sfx_N^* - 0^* )\, =\, N\tau_p\, =\, \la \sft_p ,\sfx_N^* - 0^*\rab .
\]

\begin{thm}
\label{lem:A} 
There exists a positive function $G :\bbN \to R_+ $, of an at most quadratic growth; $G (z)\leqs z^2$,  such that,
\begin{equation}
\label{eq:A}
{\rm e}^{\la \sft_p, \sfu -\sfv\rab + \la \sft_p , \sfy -\sfx \rab}
\Perc \lb \cA([\sfv ,\sfx], [\sfu ,\sfy ])\rb\, \sim\, \frac{G ( \la \sfe_2 , 
\sfx - \sfv \rab )G ( \la \sfe_2 , 
\sfy - \sfu \rab )}{N^2} , 
\end{equation}
uniformly in $\abs{\sfv}, \abs{\sfx} \leqs\log N$ and 
$\abs{\sfx_N -\sfu}, \abs{\sfx_N - \sfy} \leqs \log N$.
\end{thm}
The above function $G$ is, of course, related to renewal function $U$ 
which appears in the statement of Theorem~\ref{bigthm:B}.

\begin{lem}
\label{lem:LR}
Both sums below converge exponentially fast in $m, |\sfv | , |\sfx|$ and, respectively,
in $r, |\sfx_N -\sfu|$ 
and $|\sfx_N - \sfy |$, 
\begin{equation}
\label{eq:LR}
\begin{split}
&\sum_{m>0}\ \ \sumtwo{\sfv, \sfx\in\cH_m}{\sfv\prec\sfx}
{\rm e}^{\la \sft_p ,\sfv \rab +\la\sft_p , \sfx\rab}
\Perc\lb \cL([\sfv ,\sfx ])\rb \, \\
& =\, 
\sum_{r\geq 0}\ \ \sumtwo{\sfu, \sfy\in\cH_{N-r}}{\sfu\prec\sfy}
{\rm e}^{\la \sft_p ,\sfx_N -\sfu\rab  +\la \sft_p,\sfx_N -\sfy\rab}
\Perc\lb \cR([\sfu ,\sfy ])\rb \, <\,
\infty .
\end{split}
\end{equation}
\end{lem} 

The main effort will be to prove Theorem~\ref{lem:A}. It is precisely
at this stage we shall need the full power of the theory developed
in \cite{CI} and its geometric adjustment as in \cite{CIV3} combined
with results on asymptotic behaviour and repulsion of a pair of 
non-intersecting random walks. 
On the other hand, Lemma~\ref{lem:LR} and Lemma~\ref{lem:decomposition}
follow by  a simple adjustment of the renormalization mass-gap type bounds obtained 
in \cite{CI}. 
Accordingly, in the remaining of this subsection we shall briefly 
recall these mass-gap estimates and, subsequently, 
explain \eqref{eq:LR} and and \eqref{eq:decomposition}. The
more difficult proof of \eqref{eq:A} will be postponed to the next section.


\subsection{Structure of sub-critical connections}
\label{sub:structure1}
In this section we shall recall and reformulate the results of 
\cite{CI, CIV1, CIV2, CIV3} in a form which is convenient for 
later use. 
\subsubsection{Geometry of the inverse correlation length}
For any $p <p_c$ the inverse correlation length is defined
via
\begin{equation}
\label{Inverse_cr}
\tau_p (\sfx)\, =\, -\lim_{n\to\infty} \frac1{n}\log \Perc
\lb 0\lra \lfloor nx\rfloor \rb .
\end{equation}
As it was mentioned above the inverse correlation length at a sub-critical
$p$ equals to the surface tension at the dual super-critical 
value $p^*$. A fundamental result \cite{Me, AB} implies that
$\tau_p$ is an equivalent norm on $\bbR^2$ for every $p <p_c$. 
As such $\tau_p$ is the support function of the convex compact set
$\Kp$, which in fact is precisely the Wulff shape for the 
dual super-critical model. The relation between $\Kp$ and 
$\tau_p$ is given by
\begin{equation}
\label{Wulff}
\Kp\, =\, \bigcap_{\sfx\neq 0}
\setof{\sft\in\bbR^2}{\la\sft ,\sfx \rab\leq \tau_p (\sfx )}
\quad\text{and}\quad \tau_p (\sfx )= \max_{\sft\in\pKp} \la \sft ,\sfx \rab .
\end{equation}
Alternatively(\cite{CI}) $\Kp$ is the closure of the 
domain of convergence of the series
\begin{equation}
\label{Series}
\sft\in\text{int}\Kp\, \Longleftrightarrow\, \sum_{\sfx\in\bbZ^2} 
{\rm e}^{\la \sft ,\sfx\rab }\bbP_{p}\lb 0\lra \sfx\rb \, <\, \infty.
\end{equation}
Furthermore, as it has been proven in \cite{CI}, the boundary $\pKp$
 is locally analytic and has a strictly positive curvature. 
In particular, for each $\sfx\neq 0$ there is a uniquely defined dual 
point $\sft =\sft_{\sfx}\in\pKp$, such that
\[
\tau_p (\sfx )\,  = \, \la \sft ,\sfx \rab .
\]
Geometrically, $\sfx$ is orthogonal  to the tangent space $T_{\sft}\pKp$,
\begin{equation}
\label{Tangent}
\la\sfx , \sfv\rab =0\quad\forall\, \, \sfv\in T_{\sft}\pKp .
\end{equation}

\subsubsection{Forward cone $\cC_\delta$} 
Recall that $\sft_p = (\tau_p ,0 )\in \pKp$ is the dual point to the
horizontal axis direction $\sfe_1$.  

 Let $\delta >0$ be fixed. 
The forward cone $\cC_\delta$ is defined as follows,
\begin{equation}
\label{Deltacone}
\cC_\delta\, =\, \setof{\sfx = (t,x)\in\bbR^2}{\la \sft_p ,\sfx \rab \geq
(1-\delta)\tau_p (\sfx )} .
\end{equation}
In view of the axis symmetries and angular strict convexity 
of $\tau_p$ there exists $\alpha >0$, 
such that
\[
\cC_\delta\, =\, \setof{\sfx = (t,x)\in\bbR^2}{0\leq |x|\leq \alpha t} .
\]
It happens, however,  that the $\tau_p$-metrics naturally captures
the geometry of the problem and, accordingly, we shall stick to
the definition \eqref{Deltacone}.

\subsubsection{Cone points of $\Cl (\sfx ,\sfy )$} 
Let $\sfx ,\sfy\in \bbZ^2$  and assume that 
the cluster $\Cl (\sfx ,\sfy )\neq\emptyset$. In such a case 
 we say 
that a point $\sfz\in\Cl (\sfx ,\sfy )$  is a cone point 
of the latter  if 
$\sfz$ lies strictly between $\sfx$ and $\sfy$ with respect to the 
$\sfe_1$ direction,
\begin{equation}
\label{Conepoint1}
\sfz \in \Cl (\sfx ,\sfy )\quad\text{and}\quad
\la \sft_p , \sfx \rab < \la \sft_p , \sfz\rab < \la \sft_p , \sfy \rab ,
\end{equation}
and, in addition (Figure~{3}), 
\begin{equation}
\label{Conepoint2}
\Cl (\sfx ,\sfy )\, \subseteq\, \lb \sfz -\cC_\delta\rb
\cup \lb \sfz +\cC_\delta\rb .
\end{equation}
\begin{figure}[htb] 
\psfrag{x}{$\sfx$}
\psfrag{y}{$\sfy$}
\psfrag{z}{$\sfz$}
\psfrag{Cp}{$\sfz+\cC_\delta$}
\psfrag{Cm}{$\sfz -\cC_\delta$}
\centerline{\includegraphics{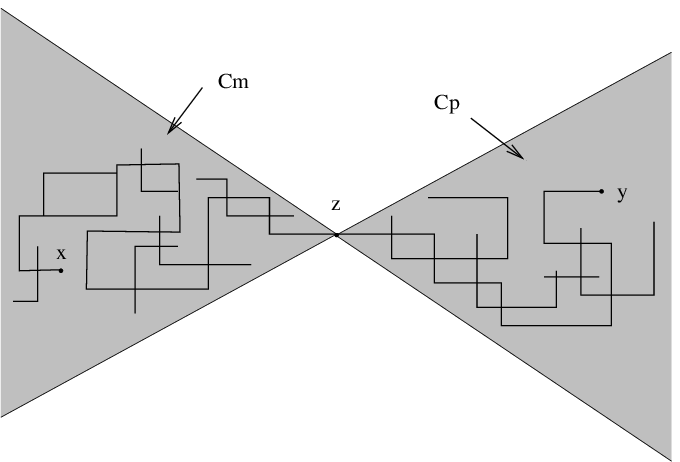}}
\label{Fig3}
\caption{$\sfz$ is a cone point of $\Cl (\sfx, \sfy )$}
\end{figure}
Clearly, $\Cl (\sfx ,\sfy )$ cannot have any cone points at all once 
$\sfy\not\in \sfx +\cC_\delta$. In the latter case, however, 
\[
\tau_p (\sfy -\sfx ) >\la t_p , \sfy -\sfx\rab +c_2\delta |\sfy -\sfx | ,
\]
where
\[
c_2 = \min_{\sfv\in\bbS^1} \frac{\tau_p (\sfv )}{ |\sfv |} .
\]
Consequently, there exists $\red{\nu_0} = \red{\nu_0}(p ,\delta )>0$ such that
\begin{equation}
 \label{eq:NotinCone}
\Perc\lb 0\lra \sfx\rb\, \leqs\, {\rm e}^{-\la \sft_p ,\sfx\rab - \red{\nu_0}\abs{\sfx}} , 
\end{equation}
uniformly in $ \sfx\not\in\cC_\delta$.

On the other hand,  for $\sfx\in\cC_\delta$, 
the techniques developed in  \cite{CI, CIV1, GI,  CIV3} readily 
imply the following mass-gap type result: For $0<k<m$ 
and $\sfx\in\cH_m$ consider the event, 
\[
 \cN_{k,m} (\sfx)\, \df\, \lbr  \Cl (0,\sfx )\ \text{has no cone points in $\cH_{k,m}$}\rbr .
\]
Then,
\begin{thm}
\label{thm:Massgap}
There exists $\red{\nu_1} = \red{\nu_1}(p ,\delta )>0$ such that uniformly in $k,l\in\bbN$ and  
in $\sfx\in\cC_\delta\cap\cH_{k+l}$, 
\begin{equation}
\label{eq:Massgap}
\Perc \lb 0 \lra\sfx ~;~\cN_{k ,k+l} (\sfx) \rb\,
 \leqs\, {\rm e}^{- \la\sft_p , \sfx \rab - \red{\nu_1} l}  .
\end{equation}
\end{thm}
{\it Proof.} A straight forward adaptation of the arguments in ~\cite{CI, CIV1, GI,  CIV3}. 
\qed

\bigskip
\noindent
Together \eqref{eq:NotinCone} and \eqref{eq:Massgap} imply: There 
exists $\red{\nu_2} = \red{\nu_2} (p ,\delta )>0$, such that uniformly in 
$l\in\bbN$ and in $\sfx\in\cH_l $, 
\begin{equation}
 \label{eq:kSumBound}
\sum_{k\geq 0}\Perc\lb 0\lra\sfx +k\sfe_1~;~\cN_{k ,k+l} (\sfx +k\sfe_1 )\rb\, 
\leqs\, {\rm e}^{-\la \sft_p ,\sfx \rab - \red{\nu_2}\abs{\red{\sfx}}} .
\end{equation}
Indeed, if $\sfx +k\sfe_1\not\in\cC_\delta$, then by \eqref{eq:NotinCone}, 
\[
 \Perc\lb 0\lra \sfx +k\sfe_1\rb \, \leqs\, {\rm exp}\lbr -\la\sft_p , \sfx +k\sfe_1\rab
 - \red{\nu_0}\abs{\sfx +k\sfe_1}\rbr \, \leq\, {\rm e}^{-\la\sft_p ,\sfx\rab - \red{\nu_0}\abs{\red{\sfx}} - k\tau_p} .
\]
If, however, $\sfx + k\sfe_1\in\cC_\delta$, then 
\[
 \la \sft_p ,k\sfe_1\rab +\red{\nu_1} l ~\geq~ \red{\nu_3} \lb k+ \abs{\sfx}\rb ,
\]
for some $\red{\nu_3} =\red{\nu_3} (p ,\delta )$, and one can rely on \eqref{eq:Massgap} in 
order to conclude that
\[
 \Perc\lb 0\lra \sfx +k\sfe_1~;~\cN_{k, k+l} (\sfx )\rb\, \leqs\, 
{\rm e}^{-\la\sft_p ,\sfx\rab - \red{\nu_3} \abs{\red{\sfx}} - k\red{\nu_3}} .
\]
It follows that, 
\begin{equation}
 \label{eq:kBound}
 \Perc\lb 0\lra \sfx +k\sfe_1~;~\cN_{k, k+l} (\sfx )\rb\, \leqs\, 
{\rm e}^{-\la\sft_p ,\sfx\rab - \min\lbr \red{\nu_0} ,\red{\nu_3}\rbr \abs{\red{\sfx}} - k\min\lbr\red{\nu_3} ,\tau_p\rbr} ,
\end{equation}
uniformly in $k ,l\in\bbN$ and in $\sfx\in\cH_l$ .
Summing over $k$ yields \ref{eq:kSumBound}.

\subsection{Proof of Lemma~\ref{lem:LR}} 
\label{sub:Prooflem:LR}
Recall that  $p_c=1/2 < p^* =1-p$ and that the 
sub-critical $p$-percolation lives on the direct lattice $\bbZ^2$ .
We claim that there exists $\red{\nu_4}=\red{\nu_4} (p) >0$ such that, 
\begin{equation}
 \label{eq:BpL}
\Perc\lb \cL([\sfv ,\sfx ])\rb\, \leqs\, 
{\rm exp}\lbr
- \la \sft_p ,\sfv + \sfx\rab - \red{\nu_4} \lb |\sfv | +|\sfx|\rb \rbr 
\end{equation}
uniformly in $l\in\bbN$ and in $\sfv ,\sfx\in\cH_l$. \eqref{eq:LR} is 
an immediate consequence.  In its turn \eqref{eq:BpL} is a mass-gap estimate
of the same type as \eqref{eq:kSumBound}. More precisely, for $k\geq 0$
define
\[
 \cL_{-k} ([\sfv ,\sfx ])\, =\, \lb
\lbr (-k ,0)\lra\sfv\rbr\circ\lbr (-k ,0)\lra\sfx\rbr\rb \cap \cL ([\sfv ,\sfx ]) .
\]
Then, by a more or less straightforward adjustment of the 
arguments leading to \eqref{eq:kBound} we infer that there exists 
$\red{\nu_4 = \nu_4 (p)}$, $\red{\nu_5} = \red{\nu_5} (p)$ such that, 
\[
\Perc\lb  \cL_{-k} ([\sfv ,\sfx ])\rb\, \leqs\, 
{\rm exp}\lbr - 2k\min\lbr \tau_p ,\red{\nu_5}\rbr   - \la \sft_p ,\sfv + \sfx\rab -
\red{\nu_4}  \lb |\sfv | +|\sfx|\rb\rbr.
\]
Since, 
\[
 \Perc\lb \cL([\sfv ,\sfx ])\rb\, \leq\, \sum_k 
\Perc\lb  \cL_{-k} ([\sfv ,\sfx ])\rb ,
\]
 \eqref{eq:BpL} follows.\qed

\subsection{Proof of Lemma~\ref{lem:decomposition}}
\label{sub:Proofdecomposition}
Lemma~\ref{lem:decomposition} follows by a very similar line of 
reasoning:
\newline
 As in the case of \eqref{eq:BpL}, 
mass-gap type estimates of \cite{CI, CIV3} imply that there exists
$\red{\nu_6}=\red{\nu_6} (p) >0$, such that
\[
\Perc\lb \cI ([\sfv , \sfx ], [\sfu ,\sfy ])\rb \, \leqs\, 
{\rm e}^{-2N\tau_p -\red{\nu_6} (|\sfv| +|\sfx| +|\sfx_N^* -\sfu |
+|\sfx_N^* -\sfy |)} \quad\text{and}\quad 
\Perc\lb\cI (\emptyset )\rb \, \leqs\, {\rm e}^{-2N(\tau_p+\red{\nu_6} )} .
\]
These are a-priori bounds: 
 Once  
Theorem ~\ref{lem:A} is established they 
 render $\Perc\lb\cI (\emptyset )\rb$
or $\Perc\lb \cI ([\sfv , \sfx ], [\sfu ,\sfy ])\rb$, with at least one 
of $|\sfv|, |\sfx| , |\sfx_N^* -\sfu |
, |\sfx_N^* -\sfy | $ being $\geqs \log N$, negligible with respect
to the right hand side of \eqref{eq:decomposition}.\qed


\section{Reduction to the Effective RW Picture}
\label{sec:Reduction}
We continue to assume that $\sfv, \sfx \in \cH_m$ and $\sfu, \sfy \in \cH_{N-r}$,
with $m < N-r$. The Lemma below explains the advantage of working with events
$\cA ([\sfv ,\sfx ], [\sfu ,\sfy ])$ and, consequently, the reasons behind 
 an introduction of modified events $\wt\cI ([\sfv ,\sfx ], [\sfu ,\sfy ])$ in 
\eqref{eq:Itilde}. 
\begin{lem}
\label{lem:toproduct}
Let $m,r,\sfv, \sfx, \sfu$ and $\sfy$ be as above. Then,  
\begin{equation}
\label{eq:toproduct}
\Perc\lb \cA ([\sfv ,\sfx], [\sfu ,\sfy ])\rb\, =\, \otimes
\Perc \lb \cA ([\sfv ,\sfx], [\sfu ,\sfy ])\rb, 
\end{equation}
where  $\otimes\Perc$ means that the clusters $\Cl_{m,N- r} (\sfx ,\sfy)$ 
and $\Cl_{m,N- r} (\sfv ,\sfu)$ are sampled {\bf independently}.
\end{lem}

\noindent
{\it Proof.} Let us decompose $\cA ([\sfv ,\sfx], [\sfu ,\sfy ])$ with respect
to realizations of $\gamma^{\rm up}(\Cl_{m,N-r} (\sfx ,\sfy ))$, 
\[
\Perc\lb \cA ([\sfv ,\sfx], [\sfu ,\sfy ])\rb\, =\, 
\sum_{\gamma}
\Perc\lb \cA ([\sfv ,\sfx], [\sfu ,\sfy ])\, ,\,
\gamma^{\rm up}(\Cl_{m,N-r} (\sfx ,\sfy ))=\gamma \rb .
\]
Using $\cA ([\sfv ,\sfx], [\sfu ,\sfy ])^\star$ for $\star = \mathbf{a}, \dots
\mathbf{e}$
to denote the events described 
\red{by conditions ${\bf a)}-{\bf e)}$ in the definition of
$\cA$} in Subsection~\ref{sub:basic}, we readily see that 
\[
\cA^{\mathbf a}\cap
\cA^{\mathbf c}\cap
\{\Cl_{m,N-r} (\sfv ,\sfu)\cap \gamma =\emptyset\}\ \text{and}\ 
\cA^{\mathbf b}\cap
\cA^{\mathbf d}\cap\{\gamma^{\rm up}(\Cl_{m,N-r} (\sfx ,\sfy ))=\gamma\} 
\]
are independent under $\Perc$.\qed

\subsection{Decomposition  of $\cA ([\sfv ,\sfx ], [\sfu ,\sfy ])$}
\label{sub:decompositionA}
In light of the previous Lemma, we may calculate probabilities using the product measure. 
Since we restrict attention to the case $m, r \leqs \log N$, 
for the sake of proving Theorem  ~\ref{lem:A} we may now assume without 
loss of generality that
$m = r =0$. Thus, 
\[
\sfv = (0,v),\ \sfx =(0, x), \ \sfu =(N,u)\ \text{and}\ \sfy =(N,y) .
\]
Given $0< l <N$  and $ \sfw, \sfz\in\cH_l $ 
let us say that $\cH_l$ is a cone cut line and, accordingly, that 
$\lbr \sfw , \sfz\rbr$ is a cone couple for 
$\lbr \Cl_{0, N} (\sfv ,\sfu ) ,\Cl_{0, N} (\sfx ,\sfy )\rbr $ 
if  $\sfw$ is a cone point 
of $\Cl_{0, N} (\sfv ,\sfu )$, whereas $\sfz$ is a cone point 
of $\Cl_{0, N} (\sfx ,\sfy )$. 

\begin{figure}[tbh] 
\psfrag{a}{$\sfx$}
\psfrag{b}{$\sfv$}
\psfrag{c}{$\sfy$}
\psfrag{d}{$\sfu$}
\psfrag{w1}{$\sfw_1$}
\psfrag{w2}{$\sfw_2$}
\psfrag{w3}{$\sfw_3$}
\psfrag{w4}{$\sfw_4$}
\psfrag{z1}{$\sfz_1$}
\psfrag{z2}{$\sfz_2$}
\psfrag{z3}{$\sfz_3$}
\psfrag{z4}{$\sfz_4$}
\psfrag{G}{$\gamma^{\rm up}(\Cl_{0,N}(\sfx ,\sfy ))$}
\psfrag{Clx}{$\Cl_{0,N} (\sfx ,\sfy )$}
\psfrag{Clv}{$\Cl_{0,N} (\sfv ,\sfu )$}
\psfrag{H0}{$\cH_0$}
\psfrag{HN}{$\cH_N$}
\psfrag{H1}{$\cH_{l_1}$}
\psfrag{H2}{$\cH_{l_2}$}
\psfrag{H3}{$\cH_{l_3}$}
\psfrag{H4}{$\cH_{l_4}$}
\centerline{\includegraphics{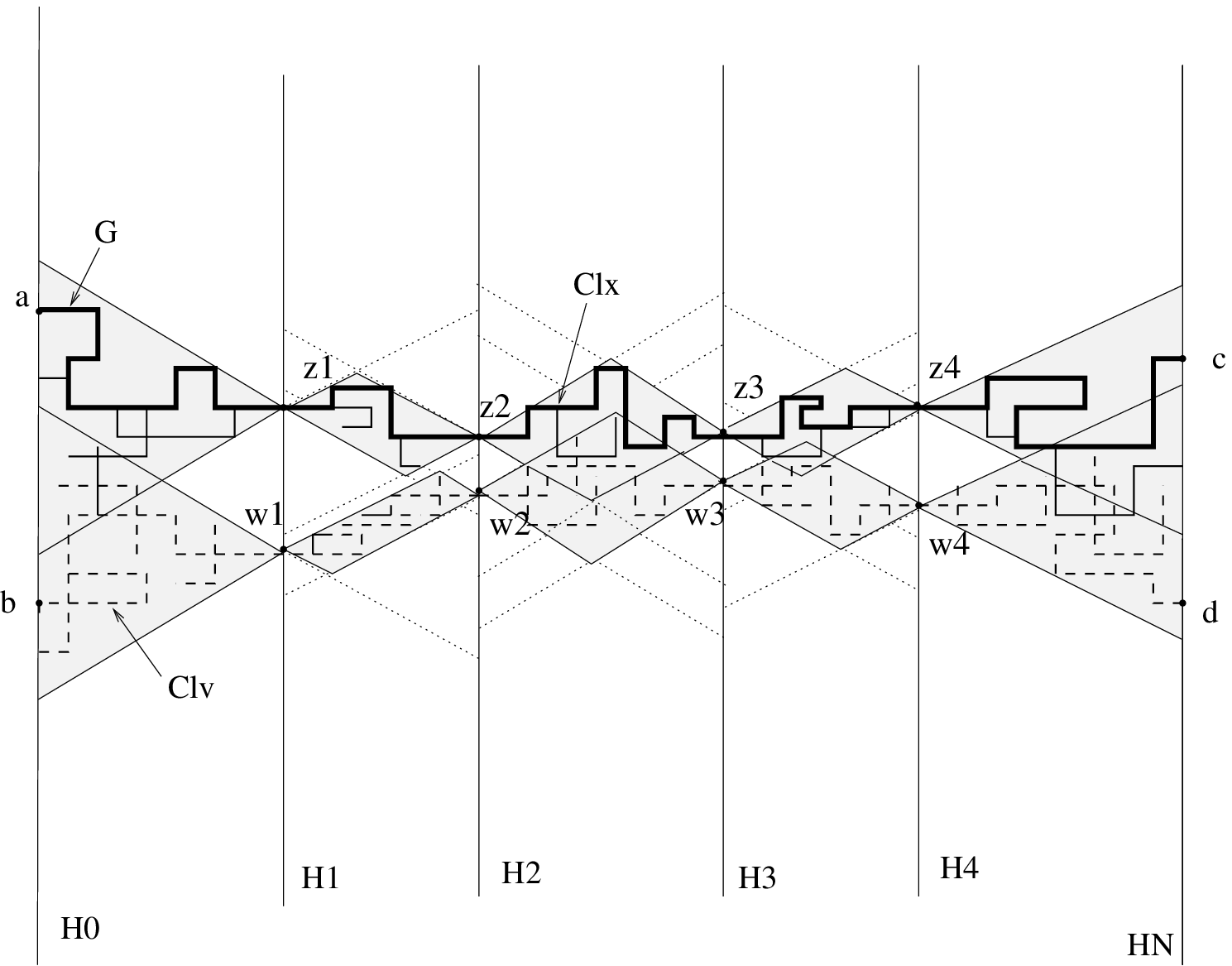}}
\label{Fig4}
\caption{Decomposition of $\cA ([\sfv ,\sfx ], [\sfu ,\sfy ])$: $l_1, l_2, 
l_3, l_4$ are cut lines.\newline
$\lbr\sfw_1 ,\sfz_1\rbr,\lbr\sfw_2 ,\sfz_2\rbr,\lbr\sfw_3 ,\sfz_3\rbr, 
\lbr\sfw_4, \sfz_4\rbr$ are 
the corresponding cone couples}
\end{figure}

A straightforward adjustment of the renormalization arguments  
behind  \eqref{eq:Massgap} in  
\cite{CI, CIV1, CIV3} implies that there exist $\red{\nu_7} = \red{\nu_7} (p, \delta ) >0$, such that, 
\begin{equation}
\label{eq:noconepoints}
\begin{split}
&\otimes\Perc\lb\lbr \Cl_{0, N} (\sfv ,\sfu ) ,\Cl_{0, N} (\sfx ,\sfy )
\rbr \, \text{has less than two
 cone cut lines}
\rb\,\\
&\qquad  \leqs\, {\rm e}^{-\la\sft_p , \sfu -\sfv\rab -\la \sft_p , \sfy -
\sfx\rab -\red{\nu_7} N} 
\end{split}
\end{equation}
uniformly in $\sfv ,\sfx ,\sfu$ and $\sfy$ under consideration.
In the case when 
$\lbr \Cl_{0, N} (\sfv ,\sfu ) ,\Cl_{0, N} (\sfx ,\sfy )\rbr$ has 
at least two cone cut lines, say $l_1 ,\dots, l_{n+1}$ with 
\[ 
\lbr \sfw_1 =(l_1 ,w_1)  ,\sfz_1 = (l_1 ,z_1) \rbr ,\dots , \lbr 
\sfw_{n+1} =(l_{n+1} ,w_{n+1} ) ,\sfz_{n+1} = (l_{n+1} ,z_{n+1}\rbr
\]
being the 
corresponding cone couples, there is a simultaneous 
irreducible decomposition (see Figure~\ref{Fig4}), 
\begin{equation}
\label{eq:irreducibleclusters}
\Cl_{0, N} (\sfv ,\sfu )\, =\, 
\Gamma_b^1\cup\Gamma_1^1\cup\dots\cup\Gamma_n^1\cup\Gamma_f^1
\quad\text{and}\quad
\Cl_{0, N} (\sfx ,\sfy )\, =\, 
\Gamma_b^2\cup\Gamma_1^2\cup\dots
\cup\Gamma_n^2\cup\Gamma_f^2 . 
\end{equation}

The simultaneous irreducible decomposition \eqref{eq:irreducibleclusters} 
sets up the stage for our effective random walk representation of 
 the double cluster 
$\lbr \Cl_{0, N} (\sfv ,\sfu ) ,\Cl_{0, N} (\sfx ,\sfy )\rbr$. In fact, our effective 
random walk will just run through the cone couples of the latter. We, therefore, 
proceed with a careful description of clusters and associated irreducible events
which show up in \eqref{eq:irreducibleclusters}

\subsection{Irreducible pairs and associated events}
\label{sub:Properties}
We shall consider the following families of clusters:

\subsubsection{Initial clusters} 
For 
 $l >0$ and $\sfw, \sfz \in \cH_l$, let $\cF_b (
[\sfw ,\sfz])$
be the set of cluster pairs $\lb\Gamma_b^1, \Gamma_b^2\rb$ satisfying:
\\ 
\noindent
(i) $\Gamma_b^i\subseteq\cH_{0,l}$ for $i=1,2$.

\noindent
(ii) $\max\lbr \Gamma_b^1\cap\cH_0\rbr = 
\max \lbr\Gamma_b^2\cap\cH_0\rbr
 =\lbr 0\rbr 
$ .

\noindent
(iii) $\Gamma_b^1\cap\cH_l =\lbr\sfw \rbr$ and $\Gamma_b^2\cap\cH_l =\lbr\sfz\rbr$. 

\noindent
(iv) $\Gamma_b^1\subseteq \sfw -\cC_\delta$ and $\Gamma_b^2\subseteq \sfz - \cC_\delta$.

\noindent
(v) $\forall k=1, \dots, l-1$, $\, \cH_k$ is not a cone 
cut line for $\lb\Gamma_b^1, \Gamma_b^2\rb$ 
    in $\cH_{0,l}$ (irreducibility).
\smallskip 

\noindent
For each such pair of clusters, with a slight abuse of notation we proceed to denote by 
$\lbr\Gamma_b^1 ,\Gamma_b^2\rbr$ the $\frE^-_{0,l}\times\frE^-_{0,l}$-measurable event that
\[
\lbr \Cl^-_{0,l} (0 ,\sfw) =\Gamma_b^1 \rbr
\times
\lbr \Cl^-_{0,l} (0 ,\sfz) =\Gamma_b^2 \rbr .
\]
Finally let:
\begin{equation}
\nonumber
\cF_b 
 = 
  \bigcup_{l \geq 1} 
    \bigcup_{\sfw, \sfz \in \cH_l} 
      \cF_b (
 [\sfw ,\sfz])
\end{equation}
In the sequel we
define random steps 
 $\sigma_b = \lb \rho_b ,\xi_b^1 ,\xi_b^2 \rb  : \cF_b\mapsto \bbZ_+\times\bbZ^2$:  If $l >0$; 
$\sfw = (l, w)$, $\sfz = (l, z)$ and $\lb\Gamma_b^1, \Gamma_b^2\rb\in\cF_b \lb [ \sfw ,\sfz ]\rb$, then
\[
 \sigma_b (\Gamma_b^1, \Gamma_b^2 ) = \lb \rho_b ,\xi_b^1 ,\xi_b^2 \rb 
= \lb l, w ,z \rb .
\]

\subsubsection{Bulk  clusters}
For 
 $l >0$ and $\sfw, \sfz \in \cH_l$, let $\cF (
[\sfw ,\sfz])$
be the set of cluster 
pairs $\lb \Gamma^1 ,\Gamma^2\rb$ satisfying:

\noindent
(i) $\Gamma^i\subseteq\cH_{0 ,l}$ for $i=1,2$.

\noindent
(ii) $\Gamma^1\cap\cH_0 =
\Gamma^2\cap\cH_0 =\lbr  0\rbr
$ 

\noindent
(iii)  $\Gamma^1\cap\cH_l =\lbr\sfw
\rbr$ and 
$\Gamma^2\cap\cH_l =\lbr\sfz
\rbr$. 

\noindent
(iv)
$\Gamma^1\subseteq \cC_\delta \cap\lb \sfw
-\cC_\delta\rb$ 
and $\Gamma^2\subseteq 
\cC_\delta
\cap\lb \sfz
-\cC_\delta\rb$ .

\noindent
(v) $\forall k=1, \dots, l-1$, $\,\cH_k$ is not a cone 
cut line for $\lb\Gamma^1, \Gamma^2\rb$ 
    in $\cH_{0,l}$ (irreducibility).
\smallskip 

\noindent
For each such pair of clusters, with a slight abuse of notation we proceed to denote by 
$\lbr\Gamma^1 ,\Gamma^2\rbr$ the $\frE^-_{0,l}\times\frE^-_{0,l}$-measurable event that   
\[
\lbr \Cl^-_{0,l} (0, \sfw
) =\Gamma^1 \rbr
\times
\lbr \Cl^-_{0,l} (0 ,\sfz
) =\Gamma^2 \rbr .
\]
Finally let:
\begin{equation}
\nonumber
\green{\cF} = 
  \bigcup_{l \geq 1} 
    \bigcup_{\sfw
, \sfz
\in \cH_l} 
      \cF ([\sfw ,\sfz]
)
\end{equation}
In the sequel we
define random steps 
 $\sigma = \lb \rho ,\xi^1 ,\xi^2 \rb  : \cF \mapsto \bbZ_+\times\bbZ^2$:  If $l >0$; 
$\sfw = (l, w)$, $\sfz = (l, z)$ and $\lb\Gamma^1, \Gamma^2\rb\in\cF \lb [ \sfw ,\sfz ]\rb$, then
\[
 \sigma (\Gamma^1, \Gamma^2 ) = \lb \rho  ,\xi^1 ,\xi^2 \rb 
= \lb l, w ,z \rb .
\]

\subsubsection{Terminal clusters}  
For 
 $l >0$ and $\sfw, \sfz \in \cH_l$, let $\cF_{\green{f}} (
[\sfw ,\sfz])$
be the set of cluster pairs $\lb\Gamma_f^1, \Gamma_f^2\rb$ satisfying: 

\noindent
(i) $\Gamma_f^i\subseteq\cH_{0, l}$ for $i=1,2$.

\noindent
(ii) $\Gamma_f^1\cap\cH_0 =
\Gamma_f^2\cap\cH_0 = \lbr 0 \rbr
$.

\noindent
(iii) $\max\lbr \Gamma_f^1\cap\cH_l\rbr =\lbr \sfw\rbr$ 
and  
$\max \lbr \Gamma_f^2\cap\cH_l\rbr =\lbr\sfz\rbr$. 

\noindent
(iv) $\Gamma_f^1\subseteq \cC_\delta\quad$ and $\Gamma_f^2\subseteq \cC_\delta$.

\noindent
(v) $\forall k=1, \dots, l-1$, $\,\cH_k$ is not a cone cut line for $\lb\Gamma_f^1, \Gamma_f^2\rb$ 
    in $\cH_{0,l}$ (irreducibility).
\smallskip

\noindent
For each such pair of clusters, with a slight abuse of notation we proceed to denote by 
$\lbr\Gamma_f^1 ,\Gamma_f^2\rbr$ the $\frE_{0,l}\times\frE_{0,l}$-measurable event that
\[ 
\lbr \Cl_{0,l} (
0 
,\sfw ) =\Gamma_f^1 \rbr
\times 
\lbr \Cl_{0,l} (0 ,\sfz )  =\Gamma_f^2 \rbr .
\]
Finally let:
\begin{equation}
\nonumber
\cF_f
= 
  \bigcup_{l \geq 1} 
    \bigcup_{\sfw, \sfz \in \cH_l} 
      \cF_f ([\sfw ,\sfz], 
)
\end{equation}
In the sequel we
define random steps 
 $\sigma_f = \lb \rho_f ,\xi^1_f ,\xi^2_f \rb  : \cF_{\green{f}} \mapsto \bbZ_+\times\bbZ^2$:  If $l >0$; 
$\sfw = (l, w)$, $\sfz = (l, z)$ and $\lb\Gamma_{\green{f}}^1, \Gamma_{\green{f}}^2\rb\in\cF_{\green{f}} \lb 
[ \sfw ,\sfz ]\rb$, then
\[
 \sigma_f (\Gamma_{\green{f}}^1, \Gamma_{\green{f}}^2 ) = \lb \rho_f  ,\xi^1_f ,\xi^2_f \rb 
= \lb l, w ,z \rb .
\]

\subsection{Construction of the effective random walk}
\label{sub:ourRW} 
Let us fix:

\noindent
1) A pair of intial clusters $\lb \wt\Gamma_b^1 ,\wt\Gamma_b^2 \rb\in\cF_b$.

\noindent
2) A sequence of pairs of clusters $\lb (\wt\Gamma_k^1 ,\wt\Gamma_k^2) \rb_{\green{{k \geq 1}}} \subset\cF$.

\noindent
3) A pair of terminal clusters $\lb \wt\Gamma_f^1 ,\wt\Gamma_f^2 \rb \in\cF_f$.
\smallskip

For $n>0$ and $\sfv = (0,v),\sfx = (0,x)$ we construct $n$-step trajectories 
of the induced effective random 
walk which starts at $\lbr \sfv ,\sfx \rbr$ as follows:  By definition, 
\begin{equation}
\label{eq:SDef}
\begin{split}
&S_0 \df (T_0 ,V_0 ,X_0) = (0, v, x),\quad 
S_k\df (T_k ,V_k ,X_k) = S_0 + \sum_1^k\sigma_l, \\
&S_k^b\df (T_k^b ,V_k^b ,X_k^b) = \sigma_b + S_k,\quad 
S_k^f\df (T_k^f,V_k^f ,X_k^f) = S_k+\sigma_f,\\[1ex]
&S_n^{bf} \df (T_n^{bf},V_n^{bf} ,X_n^{bf})  =\sigma_b +  S_n +\sigma_f,
\end{split}
\end{equation}
Above, 
\[ 
\sigma_b = \sigma_b \lb \wt\Gamma_b^1 ,\wt\Gamma_b^2\rb,\ 
 \sigma_k = \sigma \lb \wt\Gamma_k^1 ,\wt\Gamma_k^2\rb\quad \text{and}\quad  
\sigma_f= \sigma_f \lb \wt\Gamma_f^1 ,\wt\Gamma_f^2\rb .
\]
\green{Let $N>0$, $u,y \in \cH_N$ such that $S_n^{bf} = (N, u, y)$. Set also
$\sfu = (N,u)$, $\sfy = (N,y)$.}
In this notations, 
$\green{S_0}, S_0^b, S_1^b, \dots, S_n^b, S_n^{bf} = S_n^b +\sigma_f$ describes
 interpolated  trajectories through
cone cut points of  the simultaneous  irreducible decomposition of pair of clusters, 
\begin{equation}
\label{eq:EffectClusters}
\Gamma_b^1\cup\Gamma_1^1\cup\dots\cup\Gamma_n^1\cup\Gamma_f^1\quad
\text{and}\quad
\Gamma_b^2\cup\Gamma_1^2\cup\dots\cup\Gamma_n^2\cup\Gamma_f^2 ,
\end{equation}
where the induced clusters are defined as follows:
\begin{equation*}
\begin{split} 
&\lb \Gamma_b^1,\Gamma_b^2\rb = \lb \sfv +\wt\Gamma_b^1, \sfx +\wt\Gamma_b^2\rb, \\
&\lb \Gamma_k^1,\Gamma_k^2\rb = \lb (T^b_{k-1} ,V^b_{k-1}) +\wt\Gamma_k^1, \, 
 (T^b_{k-1} ,X^b_{k-1})+\wt\Gamma_k^2\rb,\\
&\lb \Gamma_f^1,\Gamma_f^2\rb = \lb (T^b_n ,V^b_n) +\wt\Gamma_f^1, \, 
 (T^b_n ,X^b_n)+\wt\Gamma_f^2\rb .
\end{split}
\end{equation*}
Comparing \eqref{eq:EffectClusters} with \eqref{eq:irreducibleclusters} we see that
 in particular 
 the above
procedure generates \green{distinctly} all the \green{cluster pairs} which contribute to the events 
$\cA\lb [\sfv ,\sfx ],[\sfu ,\sfy]\rb$ (and which have at least two cone cut lines, of course). 

Let us indroduce now the weights, 
\begin{equation*}
 \label{eq:BTilde}
\begin{split}
&\wtilde{\bbB}_p^{v,x}\lb \Gamma_b^1\cup\Gamma_1^1\cup
\dots\cup\Gamma_n^1\cup\Gamma_f^1 ; \Gamma_b^2\cup\Gamma_1^2\cup\dots\cup\Gamma_n^2\cup\Gamma_f^2\rb\\
&= \otimes\bbB_p \lb\lbr\wt\Gamma^1_b , \wt\Gamma^2_b\rbr\rb
\prod_1^n
\otimes\bbB_p \lb\lbr\wt\Gamma^1_k , \wt\Gamma^2_k\rbr\rb
\otimes\bbB_p \lb\lbr\wt\Gamma^1_f , \wt\Gamma^2_f\rbr\rb ,
\end{split}
\end{equation*}
and the events, 
\begin{equation}
\label{eq:ATilde}
\wtilde{\cA} ([\sfv ,\sfx ], [\sfu ,\sfy ]) = \bigcup_{n \green{\geq} 1} \lbr S_n^{bf} = (N,u,y) ; \quad
\Gamma^1_k \cap \gamma^{\rm up}(\Gamma^2_k) =\emptyset 
\quad \text{for} \, k=b, 1,  \dots,  n,  f \rbr .
\end{equation}
Then the irreducibility of decomposition \eqref{eq:irreducibleclusters} together with
\eqref{eq:noconepoints} imply that we can express the probability of percolation event
$\cA ([\sfv ,\sfx ], [\sfu ,\sfy ])$ under $\otimes \Perc$ asymptotically as the
probability of ``clusters-random-walk'' event $\wtilde{\cA} ([\sfv ,\sfx ], [\sfu ,\sfy ])$ 
under $\wtilde{\bbB}_p^{v,x}$:
\begin{equation}
\label{eq:BTildeOfA}
\otimes \Perc \lb \cA ([\sfv ,\sfx ], [\sfu ,\sfy ]) \rb \lb 1 +{\small O}\lb
{\rm e}^{-\red{\nu_7} N}\rb\rb = \wtilde{\bbB}_p^{v,x} \lb \wtilde{\cA} ([\sfv ,\sfx ], [\sfu ,\sfy ]) \rb .
\end{equation}
Furthermore, if we set 
\begin{equation}
\label{eq:tRplus}
\cR_n^{bf} \, \df\, \lbr X_0^b  > V_0^b \rbr \cap 
\lbr X_k^b > V_k^b : \, k=1,\, \dots,\, n \rbr \cap
\lbr X_n^{bf} > V_n^{bf} \rbr 
\end{equation}
we can write
\begin{equation}
\label{eq:ARplus}
\begin{split}
& \wtilde{\bbB}_p^{v,x}\lb\wtilde{\cA} ([\sfv ,\sfx ], [\sfu ,\sfy ])\rb \\
& =\, \sum_{n \green{\geq} 1}\wtilde{\bbB}_p^{v,x} \lb S_n^{bf} = (N,u,y)\, ;\, \cR_n^{bf}\rb\, 
\wtilde{\bbB}_p^{v,x}\lb \wtilde{\cA} ([\sfv ,\sfx ], [\sfu ,\sfy ])\, \big|\, 
 S_n^{bf} = (N,u,y)\, ;\,   \cR_n^{bf}\rb .
\end{split}
\end{equation}
We shall argue that the conditional probability above leads
 only to finite corrections, whereas sharp asymptotics are
inherited from $\wtilde{\bbB}_p^{v,x} \lb S_n^{bf} = (N,u,y)\, ;\ \cR_n^{bf}\rb$
 terms. 
This is a reduction to the effective random walk picture as
described in Subsection~\ref{sub:effective_rw}. 

\subsection{\red{Normalized step distributions}}
\label{sub:distribution}
\subsubsection{\red{Bulk} steps}
We shall now fix the steps of our effective random walk, making their 
distribution proper and check that this distribution satisifies conditions
{\bf (P1)}-{\bf (P3)} of Subsection~\ref{sub:effective_rw} whence we can
use Theorem~\ref{bigthm:B}. Let us introduce yet another probability 
measure $\bbP_{v,x}$ under which $(\sigma_k)_{k \geq 1}$ form an
infinite collection of independent random variables that share
a common distribution defined as follows:
\begin{equation}
\label{eq:P}
\bbP_{v,x}\lb \sigma = (r, x_1 , x_2)\rb\, =\, 
{\rm e}^{2r\tau_p} \wtilde{\bbB}_p^{v,x} \lb \sigma_1 = (r, x_1 , x_2)\rb 
\, = \, {\rm e}^{2r\tau_p} \sum_{\lb \Gamma^1 ,\Gamma^2\rb } \otimes\Perc\lb \lbr \Gamma^1 ,\Gamma^2\rbr \rb,
\end{equation}
where the summation is over all pairs  $\lb \Gamma^1 ,\Gamma^2\rb \in 
\cF([(r, x_1) ,(r, x_2)])$. 
We claim that $\sigma$ is a proper random variable under $\bbP_{v,x}$:  
\begin{equation}
\label{eq:Pprobability}
\sum_{(r, x_1 ,x_2 )} \bbP_{v,x}\lb \sigma = (r, x_1 , x_2)\rb\, =\, 1.
\end{equation}

Recall the notation $\sft_p = (\tau_p ,0) = \tau_p (\sfe_1 )$. Thus 
 $r\tau_p$ equals to $\la \sft_p , (r,x_i )\rab$ for $i=1,2$. For $l \leq r$
let 
us say that two points $\sfw\in \cH_l$ and $\sfz\in\cH_r$ are $c$-connected, 
$\sfw\slra{c}\sfz$ if $\Cl_{l,r}^-(\sfw ,\sfz )\neq\emptyset$ and, 
in addition, 
\[
\Cl_{l,r}^-(\sfw ,\sfz )\, \subseteq\, \lb \sfw +\cC_\delta\rb\cap
\lb \sfz -\cC_\delta\rb .
\]
The event $\lbr \sfw\slra{c}\sfz\rbr$ is $\frE_{l,r}^-$-measurable.
The results of \cite{CI, CIV3} imply the following 
consequence of \eqref{Series} : There exists a neighbourhood 
$\cU$
of $\sft_p\in \pKp$, such that for every $\sft\in\cU$, 
\[
\sum_{\sfz}{\rm e}^{\la\sft ,\sfz\rab}\, \Perc\lb 0\slra{c}\sfz\rb\, <\, 
\infty\ \Longleftrightarrow\ \sft\in\text{int}\lb\Kp\rb .
\]
As a result, for every $\sft\in\text{int}\lb\Kp\rb\cap\cU$, there
exists 
$\alpha >0$, such that
\begin{equation}
\label{eq:InKp}
\sum_{\sfz\in\cH_r}{\rm e}^{\la\sft ,\sfz\rab}\, \Perc\lb 0\slra{c}\sfz\rb\,
\leqs\,  {\rm e}^{-\alpha  r} .
\end{equation}
Conversly, for every $\sft\in\text{ext}\lb\Kp\rb\cap\cU$, there
exists
$\alpha >0$, such that
\begin{equation}
\label{eq:OutKp}
\sum_{\sfz\in\cH_r}{\rm e}^{\la\sft ,\sfz\rab}\, \Perc\lb 0\slra{c}\sfz\rb\,
\geqs\,  {\rm e}^{\alpha r} .
\end{equation}
At this point we can readily extend these convergence results to 
double clusters: 
There exists a possibly smaller neighbourhood 
$\cV\subseteq\cU$ of $\sft_p\in \pKp$, such that for every $\sft\in\cV$, 
\begin{equation}
\label{eq:phi}
\phi (\sft )\, \df\, \sum_{r>0}\, 
\sum_{\sfw ,\sfz\in\cH_r}{\rm e}^{\la \sft ,\sfw +\sfz\rab}\, \otimes\Perc
\lb 0\slra{c}\sfw , 0\slra{c}\sfz\rb\, <\, 
\infty\ \Longleftrightarrow\ \sft\in\text{int}\lb\Kp\rb .
\end{equation}
Define now, 
\[
g(\sft ) \, =\, \sum_{r>0}\, 
\sum_{\sfw ,\sfz\in\cH_r}{\rm e}^{\la \sft ,\sfw +\sfz\rab}
\sum_{\lb \Gamma^1 ,\Gamma^2\rb } \otimes\Perc\lb
\lbr \Gamma^1 ,\Gamma^2\rbr \rb ,
\]
where for each $\{\sfw ,\sfz\}$ fixed  the last 
summation is over all irreducible pairs 
$\lb \Gamma^1 ,\Gamma^2\rb \in 
\cF(\green{[\sfw  ,\sfz]})$. By 
\eqref{eq:noconepoints}
 $g$ converges everywhere
on $\cV$ (once $\cV$ is small enough).
On the other hand, for $\sft\in\cV\cap\text{int}\lb\Kp\rb$,  functions 
$\phi$ and $g$ satisfy the renewal relation, 
\[
\phi(\sft )\, =\, \frac{g ( \sft)}{1- g (\sft )} .
\]
Consequently $g (\sft )=1$ is a parametrization of $\pKp\cap\cV$.
 \eqref{eq:Pprobability} follows. \qed 

\noindent
The distribution of $\sigma$ in \eqref{eq:P}
is clearly symmetric in the sense 
of condition {\bf (P3)} of Subsection~\ref{sub:effective_rw}. 
Condition {\bf (P1)} follows from Property (iv) of
bulk clusters, as described in Subsection~\ref{sub:Properties}.
Finally, condition {\bf (P2)} is 
satisfied
by virtue of \eqref{eq:noconepoints}.

\subsubsection{Initial and terminal steps}
Analogously, we define $\sigma_b$, $\sigma_f$ to have the following distribution
under $\bbP_{v,x}$ independently of all other random variables:
\begin{eqnarray}
\nonumber
\bbP_{v,x} (\sigma_b = (r ,x_1 ,x_2 )) 
	& \df  & \bbQ_b (r, x_1 ,x_2 )\, \eqvs\, 
				{\rm e}^{2r\tau_p} \wtilde{\bbB}_p^{v,x} 
				(\sigma_b = (r ,x_1 ,x_2)) \\
\label{eq:bfsteps1}
	& = 	& {\rm e}^{2r\tau_p} 
				\sum_{\lb \Gamma_b^1 ,\Gamma_b^2\rb }\otimes\Perc \lb
				\lbr \Gamma_b^1 ,\Gamma_b^2\rbr \rb
\end{eqnarray}
and, respectively
\begin{eqnarray}
\nonumber
\bbP_{v,x} (\sigma_f = (r ,x_1 ,x_2 ))
	& \df		& \bbQ_f (r, x_1 ,x_2 )\, \eqvs\, 
					{\rm e}^{2r\tau_p} \wtilde{\bbB}_p^{v,x}
					 (\sigma_f = (r ,x_1 ,x_2)) \\
\label{eq:bfsteps2}
	& =		& {\rm e}^{2r\tau_p} 
					\sum_{\lb \Gamma_f^1 ,\Gamma_f^2\rb }\otimes
					\Perc \lb\lbr  \Gamma_f^1 ,\Gamma_f^2\rbr \rb ,
\end{eqnarray}
where the summation is over all initial irreducible pairs 
$\lb\Gamma^1_b ,\Gamma^2_b\rb\in 
\cF_b \lb [(r, x_1) ,(r, x_2)] \rb$ and, respectively, 
over all terminal irreducible pairs 
$\lb\Gamma^1_f ,\Gamma^2_f\rb\in 
\cF_f\lb  [(r, x_1) ,(r, x_2)]\rb$. 
By definition the $\eqvs$ relations in  \eqref{eq:bfsteps1} 
\eqref{eq:bfsteps2} are tuned in such a way that both 
$\bbQ_b$ and $\bbQ_f$ become probability measures. 
In addition, as it follows from \eqref{eq:noconepoints}, both display 
exponential tails. In particular, 
\begin{equation}
\label{eq:QbQfconverge}
\sum_{r}\sum_{x_1 ,x_2} \bbQ_b (r, x_1 ,x_2 )\, =\, 1\quad
\text{and}\quad 
\sum_{r}\sum_{x_1 ,x_2} \bbQ_f (r, x_1 ,x_2 )\, =\, 1
\end{equation}
converge exponentially fast in all the arguments.

\subsubsection{Trajectories}
To complete the setup, we carry over to $\bbP_{v,x}$ the definitions in \eqref{eq:SDef} and
note that under $\bbP_{v,x}$, $(S_k)_{k \geq 0}$ and $(\sigma_l)_{l \geq 1}$ satisfy the
conditions preceding Theorem~\ref{bigthm:B}. Moreover, the following holds:
\begin{equation}
\label{eq:PBConnection}
{\rm e}^{2N\tau_p} 
\wtilde{\bbB}_p^{v,x} \lb (S^b_0,\,  S^b_1,\,\dots,\, S^b_n,\, S^{bf}_n) = \underbar{s} \rb
\, \eqvs \,
\bbP_{v,x} \lb (S^b_0,\,  S^b_1,\,\dots,\, S^b_n,\, S^{bf}_n) = \underbar{s} \rb
\end{equation}
for any trajectory $\underbar{s}$ ending at time $N$.
We shall use $\bbP$   
as a short-hand notation  for $\bbP_{0,0}$.

\begin{rem}
\label{rem:Utilde}
We would like to argue that $\bbQ_b (r, x_1, x_2) = \bbQ_f (r, -x_1, -x_2)$. This does not follow immediately from symmetry with respect to reflection, since in fact, events in summation \eqref{eq:bfsteps1} are $\frE^-_{0,r}\times\frE^-_{0,r}$-measurable
while events in summation \eqref{eq:bfsteps2} are from $\frE_{0,r}\times\frE_{0,r}$ - hence not entirely symmeteric.
Nevertheless, $\otimes \Perc$-probabilities of corresponding events in the two summations differ only by a constant factor $(1-p)^4$ and thus after the $\eqvs$ normalization this difference disappears.
\end{rem}

\subsection{Proof of Theorem~\ref{lem:A}}
\label{sub:proofA}
Let us go back to \eqref{eq:ARplus}. Let $\mu = \bbE\rho$ be 
the expected value of the time coordinate displacement along an 
irreducible step (under distribution \eqref{eq:P}). First of all
note that one can restrict attention to values of
$n$ which satisfy $| N-n\mu |\leqs \sqrt{N\log N}$ .
Indeed, as it easily follows from local limit computations for $K$ sufficiently large, 
\begin{equation}
\label{eq:Krange}
\sum_{n: | N-n\mu | >K\sqrt{N\log N}}\bbP\lb \rho_1 +\dots +
\rho_n =N\rb \,=\,  
{\small o}\lb\frac1{N^2}\rb , 
\end{equation}
which is negligible with respect to the right hand side of \eqref{eq:A}.

\noindent
For $n$-s in the band $n\mu \in [N-K\sqrt{N\log N} , N+K\sqrt{N\log N}]$
and $|v|,|x|,|u|,|y| \leqs \log N$
 we proceed as follows: 

\subsubsection{Term $\wtilde{\bbB}_p^{v,x}\lb S_n^{bf} = (N, u, y); \cR_n^{bf}\rb$}
This is a purely random walk term. In view of \eqref{eq:PBConnection}  
\begin{equation}
\label{eq:Firstterm}
\begin{split}
&{\rm e}^{2\tau_p N}\wtilde{\bbB}_p^{v,x}\lb S_n^{bf} = (N, u, y); \cR_n^{bf}\rb\, 
\eqvs \bbP_{v,x}\lb S_n^{bf} = (N, u, y); \cR_n^{bf}\rb \\
&\ \ \eqvs \, 
\sum_{l, r >0}\sum_{w_1 <z_1}\sum_{w_{n+1}<z_{n+1}}
\bbQ_b (l, w_1 - v, z_1 -x)\\
&\quad\qquad\quad\quad\ 
\qquad\qquad\qquad\times \bbP_{w_1 , z_1}\lb S_n = (N-r - l, w_{n+1},
z_{n+1})\, ;\, \cR_n^+\rb \\
&\quad\qquad\qquad\qquad\qquad\qquad\qquad\ \ 
\times\bbQ_f (r ,u-w_{n+1} , y- z_{n+1}). 
\end{split}
\end{equation}
In view of the exponential tails of  $\bbQ_b$ and $\bbQ_f$ and Remark~\ref{rem:Utilde}, 
it is now a straightforward consequence of Theorem~\ref{bigthm:B} that 
\begin{equation}
 \label{RWTerm}
\begin{split}
\sum_n \bbP_{v,x}\lb S_n^{bf} = (N, u, y); \cR_n^{bf}\rb\, &\sim\, 
\frac1{N^2}
\sum_{l}\sum_{v^\prime < x^\prime}
\bbQ_b (l, v^\prime - v, x^\prime  -x)
 U(x^\prime  -v^\prime )\\
&\ \times\sum_r\sum_{u^\prime<y^\prime}
U(y^\prime -u^\prime)
\bbQ_f (r, u- u^\prime , y- y^\prime ) \\
& \df\, 
\frac{\widetilde{U}(x-v)\widetilde{U} (y-u)}{N^2}
\end{split}
\end{equation}
uniformly 
in $|v|, |x|, |u|, |y| \leqs \log N$. 

\subsubsection{Term $\wtilde{\bbB}_p^{v,x}\lb 
\wtilde{\cA} ([\sfv ,\sfx ], [\sfu ,\sfy ])\, \big|\, 
 S_n^{bf} = (N, u,y)\, ;\,   \cR_n^{bf}\rb$} We would like to argue that 
under $\cR_n^{bf}$ the trajectories of upper and lower random
walks are repulsed and, consequently, the additional constraint
$\Cl_{0, N}(\sfv , \sfu)\cap\gamma^{\rm up}(\Cl_{0, N} (\sfx ,\sfy ))$ 
imposed by the event $\wtilde{\cA} ([\sfv ,\sfx ], [\sfu ,\sfy ])$ actually
applies only close to the $\cH_0$ and $\cH_N$ lines and, furthermore, this 
constraint asymptotically decouple. 
In fact, 
 we claim:
\begin{lem}
\label{lem:Decoupling} 
There exists a positive function $H$ on $\bbZ$ of 
an at most linear growth; $H(z)\leqs z$, such that 
\begin{equation}
\label{eq:H}
\wtilde{\bbB}_p^{v,x}\lb 
\wtilde{\cA} ([\sfv ,\sfx ], [\sfu ,\sfy ])\, \big|\, 
 S_n^{bf} = (N, u,y)\, ;\,   \cR_n^{bf}\rb\, 
\sim\, H(x-v)H(y-u) 
\end{equation}
uniformly 
in $|v|, |x|, |u|, |y| \leqs \log N$ and in 
$|n\mu -N|\leqs\ \sqrt{N\log N}$. 
\end{lem}
Lemma~\ref{lem:Decoupling} is proved in the concluding Section~\ref{sec:Decoupling}.
\smallskip

\noindent
Combining \eqref{eq:H}, \eqref{RWTerm}, \eqref{eq:ARplus} and \eqref{eq:BTildeOfA}
we recover
\eqref{eq:A} with $G (\cdot )= H(\cdot )\wtilde{U} (\cdot )$. \qed


\section{Effective random walk}
\label{sec:rw}
Let $\gs_k  = (\rho_k , \xi^1_k, \xi^2_k)$ be a sequence of
i.i.d random  variables on $\bbN\times \bbZ^2$ which satisfy
conditions {\bf (P1)}-{\bf (P3)} of
Subsection~\ref{sub:effective_rw}. In the sequel we shall stick to
the notation introduced before, in particular, the event
$\cR_n^+$ is the one defined in \eqref{eq:Rplus} and 
\[
S_n = (T_n, V_n, X_n) = S_0 + \sum_{k 
 n}
\gs_k ,
\]
 is the trajectory of the random walk.  We use $\bbP_{v, x}$ for 
the distribution of the random walk with $S_0 = (0, v, x)$. 
Set,
\[
r_n (t; v,x; u, y) = \bbP_{v, x} (S_n = (t, u, y), \sas{\cR}_n) .
\]
In this notation  the left-hand side of \eqref{eq:RWRplus} equals to 
\[
 \sum_n r_n (N; v,x; u, y) .
\]
Let $p_n(t; v,x; u,y ) = \bbP_{v,x}(S_n = (t, u, y))$ be the transition
probabilities of the unconstrained  walk $S_n$.  The main computation, which 
is built upon combinatorial techniques developed in \cite{AD,BJD} (and 
is essentially contained
in those papers),
relates
$r_n$ and $p_n$: Let
$\mu = \bbE
\rho$ be the average length of a step along the time axis.

For the rest of the section fix a function $\delta :\bbN\mapsto \bbR_+$ of an 
almost linear growth,
\begin{equation}
\label{deltaGrowth}
 \forall~\alpha >0\ \lim_{n\to\infty}\frac{\delta (n)}{n^{1-\alpha }}=\infty\quad{\rm but}\quad
\lim_{n\to\infty}\frac{\delta (n)}{n}= 0.
\end{equation}
\begin{thm}
\label{bigthm:C}
Assume {\bf (P1)}-{\bf (P3)}. There exists 
a positive function $U$ on $\bbN$ of an at most linear growth;
$U(z)\leqs z$, such that for every  $\epsilon \in (0, 1/4)$,
\begin{equation}
\label{eq:AsymptoticsVarSteps}
r_n (t, v,x, u, y)\, \sim\,  \frac{U( x-v) U(y-u )}{n} \, p_n (t ; v,x ;u,y) ,
\end{equation}
uniformly in $x > v$, $y > u$ such that $\max\lbr |u-v|, |y-x|,|t -n\mu|\rbr \leqs \delta (n)$ 
and such that $\max\lbr x-v  ,y-u\rbr \leqs n^\epsilon$.
\end{thm}
Note that in the regime  $|t -n\mu| \geqs  \delta (n)$ the function 
$r_n (t, v,x, u, y)$ has an at least  stretched  exponential decay. Thereby, 
the target claim 
\eqref{eq:RWRplus} of Theorem~\ref{bigthm:B} routinely follows then from  
\eqref{eq:AsymptoticsVarSteps}, 
usual 
local limit description of $p_n$ and  Gaussian summation formula.
\begin{rem}
 There is nothing sacred in condition \eqref{deltaGrowth}. It just simplifies the 
formulas involved in the regime we actually need to apply them:
 However, since random variables $\sigma_k$ have exponential tails and since below
we shall rely {\em only} on the symmetries of 
$Z_n = X_n - V_n$ but not on the symmetries
of each of the two random walks involved, which in particular enables tiltings of the 
type $\lambda_1 T_n + \lambda_2 (X_n + V_n )$,
we could have readily extended \eqref{eq:AsymptoticsVarSteps} 
to the case of $\max\lbr |u-v|, |y-x|,|t -n\mu|\rbr < \nu n$ (for some fixed positive $\nu$)
but with, of course, appropriately modified renewal functions .
\end{rem}

We shall start by analyzing the difference $Z_n = X_n -V_n$, which is in itself 
a one dimensional random walk with symmetric steps having exponentially decaying distributions.
The event $\cR_n^+$ can be recorded in terms of $Z_n$ as 
\[
 \sas{\cR}_n\, =\, \lbr Z_k >0\ \ \text{for}\  k=1,\dots ,n\rbr .
\]
Let $\bbP_w = \bbP (~\cdot~| Z_0 =w )$, $q_n (\cdot ,\cdot )$ is the transition function
of $Z_n$, and let 
\[
 u_n (w,z )\, =\, \bbP_w \lb~ Z_n =z ; \cR_n^+ \rb . 
\]
Then,
\begin{thm}
\label{thm:ZnResult}
There exists 
a positive function $U$ on $\bbN$ of an at most linear growth;
$U(z)\leqs z$, such that for every  $\epsilon \in (0, 1/2)$:
\begin{equation}
\label{ZnResult}
u_n (w ,z )\, \sim\, \frac{U (w) U(z)}{n} q_n (w ,z ) , 
\end{equation}
uniformly in $w,z$ such that $0 < w,z \leqs n^\epsilon$. 
\end{thm}
A proof (which is based on \cite{AD,BJD}) will be given in Subsection~\ref{sub:Zn}.
The extension to Theorem~\ref{bigthm:C} will be explained in Subsection~\ref{sub:FullRW}. 
Finally, 
 Section~\ref{sec:Decoupling} 
is devoted to proofs of Proposition~\ref{prop:diamonds} 
and Lemma~\ref{lem:Decoupling}.

\subsection{One dimensional random walk $Z_n$ conditioned to stay positive}
\label{sub:Zn}
\subsubsection{Ladder  variables and Alili-Doney representation} In the 
sequel $\bbP$ is a shorthand notation for $\bbP_0$, $q_n(z)$ is 
a shorthand for $q_n(0,z)$. 

Let us say that $n$ is a strictly ascending ladder time if, 
\begin{equation}
\label{eq:Lplus}
\sas{\cL}_n = \{Z_n > Z_k \, ; \quad k=0, ...., n-1 \} 
\end{equation}
happens. 
A standard time
reversal argument (c.f. \cite{Feller.Book.Bible.Vol2}; XII, 2) implies
that under $\bbP$ the events 
$\sas{\cL}_n \cap \lbr Z_n = z \rbr$ and
$\sas{\cR}_n \cap \lbr Z_n = z \rbr$ have the same
probability for every $z>0$.

Similarly, let us say that $n$ is a non-strictly
ascending ladder time if,
\begin{equation}
\label{eq:Lnot}
\sao{\cL}_n = \lbr Z_n \geq  Z_k \, ; \quad k=0, ...., n -1\rbr.
\end{equation}
happens. Then, under $\bbP$  and for every $z\geq 0$, 
 the event $\sao{\cL}_n \cap \lbr Z_n = z \rbr$ 
 has the same
probability as the event $\sao{\cR}_n \cap \lbr Z_n = z\rbr$, where, 
\[
 \sao{\cR}_n\, =\, \lbr Z_k \geq 0\ \ \text{for}\  k=\red{0},\dots ,n\rbr .
\]
Define 
\begin{equation}
\label{eq:Nvariables}
  \begin{split} 
\sas{N} (z)\, & =\, \# \lbr m\geq 0~:~ Z_m < z\ \text{and}\ \sas{\cL}_m\rbr  \\
\sao{N}_n (z)\, & =\, \# \lbr m= 0 \dots n~:~ Z_m \leq z\ \text{and}\ \sao{\cL}_m\rbr \quad ; \quad
\sao{N} (z)\, =\, \lim_{n \to \infty} \sao{N}_n(z) . 
\end{split}
\end{equation}
The results of \cite{AD,BJD}, which are based on a beautiful generalization of Feller's 
combinatorial path surgery lemma, state:
\begin{equation}
 \label{ADResult}
\begin{split}
u_n (z)\df \bbP\lb Z_n =z ,\sas{\cR}_n\rb &= 
\bbP \lb Z_n =z ,\sas{\cL}_n\rb\, =\, \frac1{n}\bbE \lb\sas{N} (z); Z_n =z\rb\\
&\text{and}\\
\sao{u}_n (z)\df \bbP\lb Z_n =z ,\sao{\cR}_n\rb &= 
\bbP \lb  Z_n =z ,\sao{\cL}_n\rb\, =\, \frac1{n}\bbE \lb \sao{N}_n (z); Z_n =z\rb ,
\end{split}
\end{equation}
where the first identity holds for all $z>0$, whereas the second identity holds for 
every $z\geq 0$. 
\subsubsection{Apriori bounds}
Combinatorial identities \eqref{ADResult} readily yield a priori bounds on $u_n$ and $\sao{u}_n$
in terms of the unconstrained transition function $q_n$.
Indeed, since, by construction, $\sas{N}(z)\leq z$, we trivially infer:
\begin{equation}
 \label{ZSAB}
u_n (z)\, \leq\, \frac{z}{n}q_n (z) .
\end{equation}
In the case of non-strict ladder variables note that $\sao{N}(z)$ can be represented
as 
\begin{equation}
\label{Nnotsum}
 \sao{N}(z)\, =\, \sum_0^{\sas{N} (z+1)}\eta_l ,
\end{equation}
where $\eta_l$ are i.i.d. geometric random variables, independent of $\sas{N}(z+1)$,
with probability of failure 
\begin{equation}
\label{Failure} 
\chi\, =\, \bbP\lb \exists n~:~ Z_n =0\ \text{and}\ Z_m <0\ \text{for}\ m=1, \dots n-1\rb .
\end{equation}
Using H\"{o}lder inequality with $a > 1$, $a^* = 1 / (1-1/a)$, we get
\begin{equation}
  \begin{split}
    \sao{u}_n(z) & \leq \frac{1}{n} \lb \bbE \lb \sao{N}_n(z) \rb^a \rb^{1/a} 
      \lb q_n(z) \rb^{1/a^*} \\
    & \leq \frac{1}{n} (z+1) \lb \bbE \eta^a \rb^{1/a} 
      \lb q_n(z) \rb^{1/a^*} 
  \end{split}
\end{equation}
which gives for a fixed $a$
\begin{equation}
\label{ZWAB}
\sao{u}_n(z) \leqs \frac{z+1}{n} \lb q_n(z) \rb^{1/a^*} .
\end{equation}
As an a priori bound this fits in with our computations perfectly well once $a^*$ is 
sufficiently close to one. 
Using standard local limit results,  let us record \eqref{ZSAB} and \eqref{ZWAB} as
\begin{equation}
\label{ZAB}
u_n (z)\, \leqs\, \frac{z}{n^{3/2}}\quad\text{and}\quad
\sao{u}_n (z)\, \leqs\, \frac{z+1}{n^{1+1/2a^*}}.
\end{equation}

\subsubsection{Asymptotics of $u_n (z) $ and $\sao{u}_n (z)$.}
\label{subsub:asym}
It is only a short step now to derive uniform asymptotic description of $u_n$ and $\sao{u}_n$:
Let $\epsilon \in (0 ,1/2)$. We claim that uniformly in $0\leq z\leqs n^\epsilon$,
\begin{equation}
\label{Uzet}
\bbE \lb \sas{N} (z+1); Z_n =z+1\rb 
\, \nnsim{1-\chi}\, \frac{U(z+1)q_n (z+1)}{1-\chi} \, \nnsim{1}\, 
\bbE \lb \sao{N}_n (z); Z_n =z\rb ,
\end{equation}
where $U(z)$ is the renewal function
\begin{equation}
 \label{Ufunction}
U (z)\, =\, \bbE\sas{N} (z)\, =\, 
\sum_{r <z}\sum_{m} \bbP\lb Z_m =r ; \sas{\cL}_m\rb\, =\, \sum_{r<z}\sum_{m} u_m (r ) .
\end{equation}
Alternatively, in view of \eqref{Nnotsum}, the renewal function $U$ could be defined
via,
\[
\frac{1}{1-\chi}U (z+1)\, =\,   \bbE\sao{N} (z)\, =\, 
\sum_{r\leq z}\sum_{m} \bbP\lb Z_m = r; \sao{\cL}_m\rb\, =\, \sum_{r\leq z}\sum_{m} 
\sao{u}_m ( r ) .
\]
Let us prove \eqref{Uzet}. Consider first the left-most
term in \eqref{Uzet}:
\[
\bbE \lb \sas{N} (z+1); Z_n =z+1\rb\, =\, \sum_{m\leq n}\sum_{r\leq z}
\bbP\lb Z_m = r; \sas{\cL}_m\rb q_{n-m}(z-r+1).
\]
Fix $\beta \in (2\epsilon,1)$ and split the above sum into three terms with $m\in [0, n^\beta ]$, 
$m\in (n^\beta ,n- n^\beta )$ and $m\in [n- n^\beta ,n]$. 

Recall that we consider $z\leqs n^\epsilon$. Therefore, 
if $m\in [0, n^\beta ]$, then 
$q_{n-m}(z-r+1)\sim q_n (z+1)$.
Accordingly, 
\begin{equation}
\label{Region1Beg}
\begin{split}
\sum_{m\leq n^\beta }\sum_{r\leq z}
\bbP\lb Z_m = r; \sas{\cL}_m\rb q_{n-m}(z-r+1) \, &\nnsim{1}\, q_n (z+1)\sum_{m\leq n^\beta }\sum_{r\leq z}
\bbP\lb Z_m = r; \sas{\cL}_m\rb .
\end{split}
\end{equation}
Now by \eqref{ZAB}
\[
\sum_{m > n^\beta}\sum_{r\leq z} \bbP\lb Z_m = r; \sas{\cL}_m\rb \leqs (z+1)^2 n^{-\beta/2}.
\] 
Since $U(z) \sim z$ (by the Law of Large Numbers) and $\beta > 2 \epsilon$ 
it follows that 
\begin{equation}
\label{eq:sumUniCvgs}
\sum_{m \leq n^\beta}\sum_{r\leq z} \bbP\lb Z_m = r; \sas{\cL}_m\rb \nnsim{1} U(z+1)
\end{equation}
uniformly in $z \leqs n^{\epsilon}$. Hence the right term in \eqref{Region1Beg} is 
\begin{equation}
\label{Region1}
\nnsim{1} U(z+1) q_n(z+1)
\end{equation}

It remains to show that the remaining two sums are negligible. But this follows from 
our a priori bounds \eqref{ZAB} and from usual local CLT bounds 
on transition probabilities of the unconstrained random walk: 
 For $m\in (n^\beta ,n- n^\beta )$, 
\begin{equation}
\label{Region2}
\begin{split}
\sum_{m\in (n^\beta ,n-n^\beta) }\ \sum_{r\leq z}
\bbP\lb Z_m =s ; \sas{\cL}_m\rb q_{n-m}(z-r+1) \, &\leqs\, 
\sum_{m\in (n^\beta ,n-n^\beta) }
\frac{(z+1)^2}{m^{3/2} (n-m )^{1/2}} \\
&\leqs \frac{(z+1)^2}{n^{(1+\beta )/2}}.
\end{split}
\end{equation}
On the other hand, for $m\in [n- n^\beta ,n]$, 
\begin{equation}
 \label{Region3}
\begin{split}
\sum_{m\in [n-n^\beta, n] }\ \sum_{r\leq z}
\bbP\lb Z_m =s ; \sas{\cL}_m\rb q_{n-m}(z-r+1) \, &\leqs\,
\sum_{m\in [n-n^\beta, n) }
\frac{(z+1)^2}{m^{3/2} (n-m )^{1/2}}\, \\&\leqs\,
\frac{(z+1)^2}{n^{(3-\beta )/2}},
\end{split}
\end{equation}
and the right hand sides of \eqref{Region2} and \eqref{Region3}
are indeed asymptotically negligible compared to \eqref{Uzet}.

The right asymptotic relation in 
\eqref{Uzet} follows along exactly the same line of reasoning using a priori bound \eqref{ZAB}
with $a^*$ sufficiently close to $1$. \qed 

\subsubsection{Proof of Theorem~\ref{thm:ZnResult}} 
Any path $\lb Z_0,\dots ,Z_n\rb$  contributing
to $u_n (w,z)$ certainly achieves its minimal value $r = \min\lbr Z_l, l=0,\dots ,n\rbr$. 
Since $Z_n$ has a symmetric distribution 
it is enough to derive asymptotics of $u_n (w ,z)$ for $w\leq z$. 
In this case, 
$ 0 <r \leq w$.  A decomposition with respect to the first time when the minimum is hit 
leads to the following representation, 
\begin{equation}
 \label{unDecomposition}
u_n (w,z)\, =\, \sao{u}_n (z-w )\, +\, \sum_{m=0}^n\sum_{r =1}^{w-1} u_m (w-r) 
\sao{u}_{n-m} (z- r) .
\end{equation}
By \eqref{ADResult} and \eqref{Uzet}, 
\[
\sao{u}_n (z-w )\, \nnsim{1}\, \frac{ U (z-w +1)q_n (z-w)}{(1-\chi)n} .
\]
As far 
as the sum in \eqref{unDecomposition} is concerned let us fix $\beta \in (2\epsilon, 1)$ 
and consider three regimes: 
$m\in [0, n^\beta ]$, 
$m\in (n^\beta ,n- n^\beta )$ and $m\in [n- n^\beta ,n]$. 
In the middle region, 
\[
 u_m (w-r) 
\sao{u}_{n-m} (z- r)\, \leqs\, \frac{(w-r)(z-r+1)}{m^{3/2} (n-m )^{1+1/2a^*}} .
\]
As a result, the contribution of $m\in (n^\beta ,n- n^\beta )$ is 
$ \leqs\, \frac{zw^2}{n^{(3+\beta )/2a^*} }$,
which is negligible compared to \eqref{ZnResult} if $\beta > 2\epsilon$ and $a^*$ is chosen sufficiently close 
to $1$.

For $m\leq n^\beta$, we substitute 
$\sao{u}_{n-m} (z-r)\nnsim{1} U(z-r +1 )q_n (z-w )/(1-\chi )n$. 
Likewise, in the regime $m\geq n-n^\beta$ 
we substitute ${u}_m (w-r)\nnsim{1} U (w-r) q_n (z-w)/n$. Putting things together, we 
conclude (see \eqref{Ufunction}):
\begin{equation}
\label{LongComp1}
 \begin{split}
 &u_n (w ,z) \, \nnsim{1}\\
&
\frac{q_n (z-w)}{n}\lb \frac{U (z-w +1)}{1-\chi}\, +\, 
\sum_{m\leq n^\beta}\sum_{r=1}^{w-1} \lsp u_m (w-r)\frac{U (z-r+1)}{1-\chi } + \sao{u}_m (z-r)U (w-r) \rsp \rb 
\end{split}
\end{equation}
By \eqref{ZAB} and $U(z) \sim z$
\[
\sum_{m > n^\beta}\sum_{r=1}^{w-1} \lsp u_m (w-r) U (z-r+1) + \sao{u}_m (z-r) U (w-r) \rsp 
\leqs U(z) U(w) n^{\epsilon - \beta / 2a^*}
\]
Consequently, once $\beta > 2 \epsilon$, $a^*$ close to $1$, we can 
drop the constraint  $m\leq n^\beta$ in the sum on the right hand side of
\eqref{LongComp1}.
 By definition (see \eqref{Ufunction}), $\sum_m u_m (w-r ) = U (w-r +1 )- U (w-r )$ and, similarly, 
\[ 
\sum_m \sao{u}_m (z -r ) = \frac{ U (z-r +1 ) - U (z -r )}{1-\chi }.
\]
As a result we get $u_n (w ,z) \, \nnsim{1}$
\begin{equation}
\label{LongComp2}
\begin{split}
& \frac{q_n (z-w)}{(1-\chi )n} \lbr \rule[-.4cm]{0cm}{.8cm} \right .
U (z-w +1)  \\
&+
\sum_{r=1}^{w-1} [U (w-r+1 ) - U(w- r)]U (z- r+1) + [U (z-r+1  ) - U(z-r)]
U(w-r ) \left . \rule[-.4cm]{0cm}{.8cm} \rbr \\
&=\, \frac{q_n (z-w)}{(1-\chi )n}
\lbr U (z-w +1)
 + \sum_{r=1}^{w-1}\lb U (w-r+1 )U (z- r+1) - U(z-r )U(w-r)\rb
\rbr
\\
&=\, \frac{q_n (z-w)}{(1-\chi )n} U(w) U(z) ,
\end{split}
\end{equation}
where on the last step we have used an  obvious relation $U (1) = 1$.\qed

\begin{rem}
\label{rem:AprioriBdTwoVars}
Theorem~\ref{thm:ZnResult} yields sharp asymptotics whenever $0\leq w,z\leq n^\epsilon$.
By using aprioi bounds \eqref{ZAB}, one can eaisly obtain from \eqref{unDecomposition} the 
following a priori bound on $u_n(w,z)$ 
 which holds uniformly in $n$, $w>0$, $z>0$: Fix 
  $a^*$ to be sufficiently close to $1$. Then, 
\begin{equation}
\label{eqn:AprioriBdTwoVars}
u_n(w,z) \leqs \frac{wz \min\{w,z\}}{n^{1+1/2a^*}} .
\end{equation}
\end{rem}

\subsection{Adjustments for $S_n$}
\label{sub:FullRW}
Let us return to the coupled RW $S_n =(T_n ,V_n ,X_n)$.  Recall that 
$Z_n = X_n-V_n$.  As in the previous subsection the events 
$\sas{\cR}_n$ and $\sao{\cR}_n$ are formulated in terms of $Z_n$. 
As usual $\bbP$ stands for $\bbP_{0,0}$ 

\subsubsection{Alili-Doney representation} Since the representation of 
\cite{AD} is based on a combinatorial identity related to a surgery of 
$Z_n$-paths, this part has an immediate generalization to the 
full $S_n$-case:
\begin{equation}
 \label{SADResult}
\begin{split}
& r_n (t; u,y)\, \df\, \bbP\lb S_n = (t; u, y); \sas{\cR}_n\rb\, =\, \frac1{n}
\bbE\lb \sas{N} (y-u); S_n = (t; u, y)\rb \quad ; \quad y > u \\
&\qquad {\rm and}\\
& \sao{r}_n (t; u,y)\, \df\, \bbP\lb S_n = (t; u, y); \sao{\cR}_n\rb\, =\, \frac1{n}
\bbE\lb \sao{N}_n (y-u); S_n = (t; u, y)\rb \quad ; \quad y \geq u .
\end{split} 
\end{equation} 

\subsubsection{Apriori bounds} In place of \eqref{ZSAB}, \eqref{ZWAB}
we now have
\begin{equation}
\label{SSAB}
r_n (t; u,y)\, \leq\, \frac{y -u}{n} p_n (t; u ,y ),
\end{equation}
and
\begin{equation}
 \label{SWAB}
\sao{r}_n (t; u,y)\, \leqs\, \frac{y -u+1}{n} \lb p_n (t; u ,y )\rb^{1/a^*},
\end{equation}
where we use a shorthand notation $p_n (t; u ,y ) = p_n (t;0,0;u,y )$ and 
$a^*$ is a fixed number as close to $1$ as needed. This follows by identical
arguments.

\subsubsection{Asymptotics of  $r_n (t; u,y)$ and $\sao{r}_n (t; u,y)$}
Fix 
 $\epsilon \in (0,1/4)$.
\red{We shall prove:}
\begin{equation}
 \label{rnBounds}
r_n (t; u,y )\, \nnsim{1}\, \frac{U (y - u )}{n}p_n (t; u, y ) \quad{\rm and}\quad
\sao{r}_n (t; u,y )\, \nnsim{1}\, \frac{U (y -u+1 )}{n (1-\chi )}p_n (t; u, y ) ,
\end{equation}
uniformly in
\begin{eqnarray}
\label{eqn:rnBoundsCond1}
& & |u|, |y|, |t-n\mu | \leqs \delta (n) \quad \text{and} \\ 
\label{eqn:rnBoundsCond2}
& & |u-y | \leqs n^{\epsilon }.
\end{eqnarray}
\red{
However, let us first assume, in place of \eqref{eqn:rnBoundsCond1} the stronger
condition that 
\begin{equation}
\label{eqn:rnBoundsCond1A}
\bbE T_n = n\mu =t \quad \text{and} 
	\quad \bbE V_n =\bbE X_n =u.   
\end{equation}
To permit the latter, we no longer suppose axes-symmetry for the distribution of
$(\xi^1 ,\xi^2)$ as required by property {\bf (P3)}. We still, nonetheless, assume
diagonal symmetry and of course {\bf (P1)} and {\bf (P2)}.
}

Set $z = y-u$. Then, \red{starting with $r_n$}, 
\begin{equation}
\label{Npexpansion}
 \begin{split}
 &\bbE\lb \sas{N} (y-u ); S_n =(t; u, y)\rb\\
&\ \sum_{m\leq n}\sum_{r =1}^{z-1}\sum_{x_m-v_m =r}\sum_s
\bbP\lb S_m = (s; v_m ,x_m ); \sas{\cL}_m\rb
p_{n-m}(t-s ; v_m,  x_m ; u, y ) \\
&=\, \sum_{m\leq n}\sum_{r =1}^{z-1}\sum_{x_1-v_1 =r}\sum_s
r_m (s; v_m ,x_m )p_{n-m}(t-s ; v_m, x_m ; u, y ),
 \end{split}
\end{equation}
where  ladder event $\sas{\cL}_m$ are still defined in terms of $Z$-process.

Now, if $\bbE T_n = t$ ,  $\bbE V_n =\bbE X_n =u$  and $|y-u| \leqs n^\epsilon$ hold 
($\epsilon < 1/4$), then 
\begin{equation}
 \label{A}
 p_n (t; u, y)\sim n^{-3/2}.
\end{equation}
As in the one-dimensional case we shall
split the sum over $m$
into three terms according to $m\leq n^\beta$, $m\in (n^\beta , n-n^\beta)$ 
and $n-m\leq n^\beta$ with $\beta \in (2\epsilon, 1/2)$.

\noindent
$\bullet$ In the region $m\leq n^\beta$ we may restrict attention to 
 $|x_m |,|v_m |,s\leqs n^\beta$. Since we choose 
$\beta <1/2$, 
\begin{equation}
\label{B} 
p_{n-m} \lb t-s ; v_m ,x_m ;u, y\rb \, \nnsim{1}\, 
p_n (t; u, y )
\end{equation}
uniformly in the remaining range of parameters. 
Hence, the corresponding contribution to the right hand-side of \eqref{Npexpansion} is,
\begin{equation}
 \begin{split}
\label{C}
&\nnsim{1} p_n (t; u, y )\sum_{m\leq  n^\beta}\sum_{r =1}^{z-1}
\sum_{x_m -v_m =r}\sum_s\bbP\lb S_m =(s ; v_m ,x_m ); \sas{\cL}_m\rb\\
&\nnsim{1}p_n (t; u, y )\sum_{m\leq  n^\beta}\sum_{r =1}^{z-1} \bbP\lb Z_m =r ; \sas{\cL}_m\rb\, \nnsim{1}\, p_n (t; u, y ) U(z ).
 \end{split}
\end{equation}

\noindent
$\bullet$ In the region $n^\beta <m < n-n^\beta $ consider \eqref{SSAB}, 
\begin{equation*}
\begin{split}
&\sum_{v_m}\sum_s
r_m (s; v_m , v_m +r )p_{n-m}(t-s ; v_m, v_m +r ; u, y )\,\\
&\quad  \leqs\, 
\frac{r}{m} 
 \sum_{v_m }\sum_s 
 p_m (s; v_m ,v_m +r )p_{n-m}(t-s ; v_m, v_m+r  ; u, y )
\end{split}
\end{equation*}
Define $\mu_{\sf T} = t/n = \bbE\rho$ and $\mu_{\sf X} = u/n = \bbE\xi^i$. Set 
$\phi_l (x) = \min\lbr |x| ,x^2 /l\rbr$ 
Since $S_l$ obeys classical local limit description
under Cramer's condition, there exists $\nu >0$, such that, 
\begin{equation}
\label{D} 
p_l (a ;b ,c)\, \leqs\, \frac1{l^{3/2}} {\rm exp}\lbr -\nu \lb \phi_l (a - l\mu_{\sf T} )
+\phi_l (b - l\mu_{\sf X} ) +\phi_l (c - l\mu_{\sf X} )\rb\rbr , 
\end{equation}
uniformly in $l, a, b$ and $c$. Consequently, 
\begin{equation}
\label{E}
 \sum_{v_m }\sum_s 
p_m (s; v_m ,r + v_m)p_{n-m}(t-s ; v_m, r +v_m ; u, y )\, \leqs\, 
\frac{1}{m^{3/2} (n-m)^{3/2}}\, \min\lbr m ,n-m \rbr, 
\end{equation}
as it follows from Gaussian summation formula.  Accordingly, the contribution to
\eqref{Npexpansion} which comes from the region $n^\beta < m< n-n^\beta$ is, 
\begin{equation}
\label{F}
\leqs\,  z^2\, \sum_{m=n^\beta}^{n-n^\beta}
\frac{\min\lbr m ,n-m \rbr}{m^{5/2} (n-m)^{3/2}}\, \sim\, 
\frac{z^2}{n^{\beta /2}}\frac1{n^{3/2}}\, \sim\, 
\frac{z^2}{n^{\beta /2}} p_n (t ;u, y).
\end{equation}
%
Since  
$z\leqs U (z)$ , 
the latter 
expression is negligible with respect to $ U (z) p_n (t ;u, y)$ as soon as 
$z \leqs n^\epsilon \ll n^{\beta /2}$. This explains the restrictions on $\beta$.
\smallskip 

\noindent
$\bullet$ In the region $n-m\leq n^\beta$ we are entitled to restrict attention to 
$|u- v_m |, |y- x_m |, t-s \leqs n^\beta$.  In such a case, $p_m (s ; v_m ,x_m )\nnsim{1}
p_n (t ;u ,y)$. On the other hand, 
\begin{equation}
\label{G}
 \sum_{v_m}\sum_s p_{n-m} (t-s ; v_m ,v_m +r; u ,v)\, \leqs\, \frac{1}{(n-m +1)^{1/2}}.
\end{equation}
Consequently the corresponding contribution to \eqref{Npexpansion} is 
$\leqs z^2 p_n(t ;u ,y) n^{\beta/2 -1}  $, which is negligible as soon as $\epsilon +\beta/2 <1$,
which is  the case if $\beta > 2 \epsilon$
\smallskip

\noindent
The $\sao{r}_n$-case could be worked out in a completely similar fashion once
 $a^*$ in \eqref{SWAB} is chosen to be sufficiently close to one.  
\red{
This proves \eqref{rnBounds} with condition
\eqref{eqn:rnBoundsCond1A} in place of \eqref{eqn:rnBoundsCond1}
}.

\subsubsection{Tilts by $\lambda = (\lambda_{\sf T}, \lambda_{\sf V}, \lambda_{\sf V})$}
\red{ 
We no longer assume \eqref{eqn:rnBoundsCond1A}, but rather
\eqref{eqn:rnBoundsCond1}, \eqref{eqn:rnBoundsCond2} and, of course, {\bf (P1)}-{\bf (P3)}.
Given $n,u ,y$ and $t$ satisfying \eqref{eqn:rnBoundsCond1}, \eqref{eqn:rnBoundsCond2}},
let us tilt $\sigma$
%
by 
an appropriate $\lambda = \lambda (n ,t, u)  = (\lambda_{\sf T}, \lambda_{\sf V}, \lambda_{\sf X})$ with 
$\lambda_V = \lambda_X$, such 
that the tilted distribution $\bbP_\lambda$ of $\sigma = (\rho ,\xi^1 ,\xi^2)$:
\be
\label{Tilt}
 \bbP_\lambda\lb \sigma = (a,b,c )\rb\, \df\, \frac{{\rm e}^{\lambda_{\sf T}a +
\lambda_{\sf X} (b+c)}}{\bbE {\rm e}^{\langle \lambda ,\sigma \rangle}}
\bbP \lb \sigma = (a,b,c )\rb
\ee
satisfies $\bbE_\lambda \sigma = (t/n ,u/n ,u/n )$. Note that in view of the symmetries
of the original $\bbP$, exponential tails of $\sigma$ and in view of \eqref{deltaGrowth}
 such tilting is always possible \red{and} as $n\to\infty$, $|\lambda | = \so$ uniformly
in the range of the parameters involved.
   
\red{
On the other hand, under $\bbP_{\gl}$ for any $\gl = (\gl_{\sf T}, \gl_{\sf V}, \gl_{\sf V})$
close enough to zero, the distribution of $\gs$ satisfies properties {\bf (P1)}, {\bf (P2)} and
the diagonal symmetry in property {\bf (P3)}. Consequently,
if we let $t_{\gl} = n \bbE_{\gl} \rho$, $u_{\gl} = n \bbE_{\gl} \xi^1 = n \bbE_{\gl} \xi^2$:}
\begin{eqnarray}
\label{eqn:rnglBounds1}
r_{n ,\lambda } (t_\lambda ; u_\lambda ,y ) 
	& \nnsim{1} 	& \frac{U_\lambda (y - u_\lambda  )}{n} p_{n ,\lambda} (t_\lambda ; u_\lambda , y ), \\
\label{eqn:rnglBounds2}
\sao{r}_{n ,\lambda} (t_\lambda ; u_\lambda ,y ) 
	& \nnsim{1} 	& \frac{U_\lambda  (y -u_\lambda +1 )}{n (1-\chi_\lambda  )}p_{n,\lambda } 
		(t_\lambda ; u_\lambda , y ) 
\end{eqnarray}
uniformly in $0 \leq y -u_{\gl} \leqs n^{\epsilon}$
\red{
with $r_{n, \gl}$, $r_{n, \gl}^0$, $U_{\gl}$, $\chi_{\gl}$
defined as in \eqref{SADResult}, \eqref{Ufunction}, \eqref{Failure}, but with
$\bbP_{\gl}$, $\bbE_{\gl}$ in place of $\bbP$, $\bbE$.
}

Furthermore, if we fix $\kappa >0$ sufficiently small.  
then the bounds \eqref{SSAB}, \eqref{SWAB} and \eqref{A}--\eqref{G} (with
\red{$u=u_{\gl}$, $t=t_{\gl}$}) also hold uniformly
for the whole family of tilted measures $\lbr \bbP_\lambda\rbr_{|\lambda |\leq \kappa}$.
\red{
Therefore, we infer that \eqref{eqn:rnglBounds1}, \eqref{eqn:rnglBounds2} also hold
}
uniformly in $|\lambda |\leq \kappa$. 

Since, 
in addition, 
\[
 \frac{p_{n,\lambda} (t; u, y)}{r_{n,\lambda} (t; u, y)} \equiv 
 \frac{p_{n } (t; u, y)}{r_{n} (t; u, y)}\quad\text{and}\quad
 \frac{p_{n,\lambda} (t; u, y)}{\sao{r}_{n,\lambda} (t; u, y)} \equiv 
 \frac{p_{n } (t; u, y)}{\sao{r}_{n} (t; u, y)} ,
\]
it suffices to check:
 %
\begin{prop}
\label{prop:ContAndUnif}
$ $

\begin{enumerate}
\item	As $\gl \to 0$, $\chi_{\gl} \to \chi_0 = \chi$
\item	As $\gl \to 0$, $U_{\gl}(z) \nnsim{1}  U_0(z) = U(z)$ uniformly in $z > 0$.
\end{enumerate}
\end{prop}
{\it Proof.}
To avoid ambiguities let us fix $\kappa >$ small enough and consider 
$\lambda = (\lambda_T ,\lambda_V ,\lambda_V )$ with $|\lambda |\df \sqrt{\lambda_T^2 +
\lambda_V^2}\leq \kappa$.  For such $\lambda$ we define tilted distributions 
$\bbP_\lambda$ as in \eqref{Tilt}. 

For $\chi_{\gl}$, write:
\[
\chi_{\gl} = \sum_{n \geq 1} \bbP_{\gl} \lb Z_n=0,\ Z_k < 0 ; k=1, \dots, n-1 \rb .
\]
For each $m$ fixed,
\be
\label{chilBound}
\begin{split} 
&\sum_{n =1}^m \bbP_{\gl} \lb Z_n=0,\ Z_k < 0 ; k=1, \dots, n-1 \rb \leq 
\chi_\lambda \\
&\quad\leq  \sum_{n =1}^m \bbP_{\gl} \lb Z_n=0,\ Z_k < 0 ; k=1, \dots, n-1 \rb +
\bbP_\lambda (\cR_m^+ ) .
\end{split}
\ee

By \red{\eqref{ZSAB}} 
\[
 \bbP_\lambda (\cR_m^+ ) \leq \frac{1}{m}\bbE_\lambda\max\lbr Z_m ,0\rbr ,
\]
which is $\leqs 1/\sqrt{m}$ uniformly in $m$ and $|\lambda |\leq \kappa$. 
 On the other hand, for every $n$ fixed the map 
$\lambda\mapsto \bbP_{\gl} \lb Z_n=0,\ Z_k < 0 ; k=1, \dots, n-1 \rb$ is evidently 
continuous. This proves (1).

In order to prove (2) consider the following probability distribution on $\bbZ_+$, 
\[
 f_\lambda (r) = \bbP_\lambda \lb \exists~m~:~ Z_1, \dots ,Z_{m-1}\leq 0\ \text{and}\  
Z_m =r\rb .
\]
 The renewal function $U_\lambda$ is recovered from $f_\lambda$ in the following way:
 Define $u_\lambda (0) =1$ and 
\be
\label{rlambda}
u_\lambda (z) = \sum_{r=1}^z f_\lambda (r) u_\lambda (z-r ).
\ee
Then, $U_\lambda (z) = \sum_0^{z-1} u_\lambda (z)$.

We claim that if $\kappa$ is sufficiently small then the family of distributions 
$\lbr f_\lambda\rbr_{|\lambda |\leq \kappa}$ has uniform exponential tails. Indeed, 
since under $\lbr \bbP_\lambda\rbr_{|\lambda |\leq \kappa}$ the distribution of steps $Z_i$ has uniform exponential tails there exists  $c_1 >0$ such that
\[
 f_\lambda (r) \leq {\rm e}^{-c_1 \red{r}} \sum_{n=1}^{\red{m}} \bbP_\lambda\lb \sao{\cR}_{n\red{-1}}\rb +
\bbP_\lambda\lb \sao{R}_m \rb \leqs {\rm e}^{-c_1 \red{r}}\sqrt{m} +\frac{1}{\sqrt{m}} .
\]
It remains to take $m = m (r) = {\rm e}^{c_1 r}$.

Standard Renewal Theory reinforced with such uniform exponential decay 
implies that as $r \to \infty$, 
\begin{equation}
\label{eqn:Limitu}
u_{\gl}(r) \to \frac1{\mu^+_{\gl}} 
\end{equation}
uniformly exponentially fast (on $|\lambda |\leq \kappa$) , 
where $\mu^+_{\gl}$ is the expected value of the strict ladder height associated with
$Z_n$, namely:
$
\mu^+_{\gl} = 
\sum_r r f_\lambda (r) .
$ .
Since, uniform exponential tails of $\lbr f_\lambda\rbr$ 
\red{and continuity of $\gl \mapsto f_{\gl}(r)$ for all $r$}  
imply that 
$\mu^+_{\gl}$ is continuous on  $|\lambda |\leq \kappa$  and since 
$\lambda\mapsto u_\lambda (z)$ is trivially continuous for every $z$ fixed, (2) follows. 
\qed

\smallskip
\subsubsection{Proof of \eqref{eq:AsymptoticsVarSteps}} 
It is enough to consider 
only the case 
of $0<w\df x-v \leq y-u \df z$. Decomposing with repect to the position of 
the first global minimum of $Z_n$, we arrive to the following 
generalization of \eqref{unDecomposition}, 
\begin{equation*}
 \begin{split}
  r_n (t; v,x;u,y )\, &=\, \sao{r}_n (t; u-v ,y-x)\\
&+\, \sum_{m=1}^{n-1}\sum_{r =1}^{w-1}\sum_{s=1}^{t-1}
\sum_{x_m -v_m =r}\hat{r}_m (s; v-v_m , x-x_m )\sao{r}_{n-m} (t-s ; u-v_m ,y -x_m) ,
 \end{split} 
\end{equation*}
where $\hat{r}_m$ is defined exactly as $r_m$ but for the reversed walk
$\hat{S}_m$ with i.i.d. steps, 
\begin{equation}
\label{sigmahatRW}
 \hat{\sigma}_k\, =\, (\rho_k ,-\xi_k^1 ,-\xi_k^2 ) .
\end{equation}
Note, however, that the distribution of $Z_n =X_n -V_n$ is always symmetric.
In particular $\hat{Z}_n$ has the same renewal function $U$ as $Z_n$. 

As before, we, applying if necessary appropriate tilts $\lambda$, may
assume that $(t ,u-v) = \bbE (T_n, V_n)$
or, equivalently, that $(t ,v-u)= \bbE (T_n, \hat{V}_n)$.
Fix $2\epsilon < \beta < 1/2$ and split the  sum over $m$ 
into three regions $m\leq n^\beta$, $n^\beta <m < n-n^\beta$ and $n-m\leq n^\beta$.
\smallskip

\noindent
$\bullet$ In the region $m\leq n^\beta$  we can restrict attention to
 $|v_m -v|, |x_m -x|$ and $s$ all being
$\leqs n^\beta$. Then, the 
second of \eqref{rnBounds} implies that
\[
 \sao{r}_{n-m} (t-s , u-v_m ,y-x_m )\, \nnsim{1}\, 
\frac{U(z-r +1)}{n (1-\chi)}p_n (t; u-v ,y-x )
\]
uniformly in $r =1, \dots ,w-1$ and such $s, v_m$ and $x_m$ with $x_m -v_m =r$.
On the other hand, for every $K >0$ fixed,
\[
\sum_{s\leq K n^\beta}\ \sumtwo{|v_m -v|, |x_m -x | \leq n^\beta}{x_m -v_m =r} 
\hat{r}_m (s; v-v_m , x-x_m )\, \nnsim{1}\, \bbP\lb 
\hat{Z}_m = w-r ; \sas{\cL}_m\rb\, =\, u_m (w -r ) .
\]
$\bullet$  
Similarly, for $n-m\leq n^\beta$ we may restrict attention to 
$|u-v_m|, |y-x_m|$ and $t-s$ being 
$\leqs n^\beta$ and, accordingly, infer 
from the first of \eqref{rnBounds} that, 
\[
 \hat{r}_{m} (s , v-v_m ,x-x_m )\, \nnsim{1}\, 
\frac{U(w-r )}{n }\hat{p}_n (t; v-u ,x-y )\, =\, 
\frac{U(w-r )}{n }{p}_n (t; u- v ,y-x ), 
\]
whereas, 
\begin{equation*}
\begin{split} 
&\sum_{t-s \leq K n^\beta}\ \sumtwo{|v_m-u|, |x_m-y | \leq n^\beta}{x_m -v_m =r} 
\sao{r}_{n-m} (t-s; u-v_m , y-x_m )\, \\
&\qquad \nnsim{1}\, \bbP\lb Z_m = z-r ; \sao{\cL}_m\rb\, =\, 
\sao{u}_{n-m} (z -r ) .
\end{split}
\end{equation*}
$\bullet$  
 As in the one-dimensional case, 
a priori bounds \eqref{SSAB} and \eqref{SWAB} (applied with $a^*$ being 
sufficiently close to one) render the contribution of the middle sum negligible. 
\smallskip

\noindent
The rest of the proof is a repetition of \eqref{LongComp1}, 
\eqref{LongComp2} and 
\eqref{eq:AsymptoticsVarSteps} follows.\qed

\subsection{Boundary steps and semi-infinite walks}
\label{sub:semi}
Assume now that $\gs_b, \gs_f$ are defined as well and have distributions
$\bbQ_b$ \eqref{eq:bfsteps1} and $\bbQ_f$ \eqref{eq:bfsteps2} under $\bbP_{v,x}$. Since the distribution $\bbQ_b$ 
of the initial step and, respectively, the 
distribution $\bbQ_f$ of the final step have exponentially decaying tails
it is straightforward to incorporate them into Theorem~\ref{bigthm:C}. 
With the random walk notation of Subsection~\ref{sub:ourRW}, set:
\[
 \tilde{p}_n (t ; v,x; u,y)\,  =\, \bbP_{v,x}\lb S_n^{bf} = (t ,u, y) \rb,
\]
and, accordingly,
\[
\tilde{r}_n (t ; v,x; u,y)\, =\, \bbP_{v,x}\lb S_n^{bf} = (t ,u, y) ;\cR^{bf}_n \rb . 
\]
Then, by \eqref{eq:AsymptoticsVarSteps} and by the very same computation as 
in \eqref{RWTerm}, 
\begin{equation}
 \label{rntilda}
\begin{split}
\tilde{r}_n (t ; v,x; u,y)\, &\sim\, 
\frac{\wt{U }(x-y)\wt{U} (y-u )}{n}
 ~\tilde{p}_n (t ; v,x; u,y) , 
\end{split}
\end{equation}
uniformly in $x > v$, $y > u$ such that $\max\lbr |u-v|, |y-x|,|t -n\mu|\rbr \leqs \delta (n)$ 
and such that $\max\lbr x-v  ,y-u\rbr \leqs n^\epsilon$,
with $\widetilde{U}$ defined precisely as in \eqref{RWTerm} and $\gep \in (0, 1/4)$.   

\noindent
Below we shall also need asymptotics of coupled random walks which take into 
account only one of $\sigma_b$ or $\sigma_f$ boundary steps. To this end let us
introduce the following notation:
\[
\tilde{p}^b_n  (t ; v,x; u,y)\, =\, \bbP_{v,x}\lb (S_n^b = (t,u ,y)\rb\, =\, 
\sum_{s ,v_b ,x_b }\bbQ_b (s ,v_b -v ,x_b - x)p_n (t-s ,u- v_b ,y -x_b ) .
\]
Similarly,  define, 
\[
 \tilde{p}^f_n  (t ; v,x; u,y)\, =\, \bbP_{v,x}\lb (S_n^f= (t,u ,y)\rb\, =\, 
\sum_{s ,v_n,x_n }p_n (s ,v_n- v ,x_n -x ) \bbQ_f (t-s ,u - u_n , y- x_n).
\]
\noindent
The corresponding versions of path non-intersection events are, 
\begin{equation}
\begin{split}
\cR_n^b \, & \df\, \lbr X_0^b  > V_0^b \rbr \cap 
\lbr X_k^b > V_k^b : \, k=1,\, \dots,\, n \rbr \\
\cR_n^f \, & \df\, \lbr X_k > V_k : \, k=1,\, \dots,\, n \rbr \cap
\lbr X_n^{f} > V_n^{f} \rbr 
\end{split}
\end{equation}

Set $\tilde{r}_n^b (t ;v,x ;u,y ) = \bbP_{v,x}\lb S_n^b = (t, u ,y) ; \cR_n^b\rb$
 and, accordingly, 
\[ 
\tilde{r}_n^f(t ;v,x ;u,y ) = \bbP_{v,x}\lb S_n^f = (t, u ,y) ; \cR_n^f\rb .
\]
Then, exactly as in \eqref{rntilda} above, 
\begin{equation}
 \label{rntildebf}
\begin{split}
 \tilde{r}_n^b (t ;v,x ;u,y )\, &\sim\, \frac{\wt{U} (x-v ) {U} (y-u)}{n} ~
\tilde{p}_n^b (t ;v,x ;u,y )\\
&{\rm and}\\
\tilde{r}_n^f (t ;v,x ;u,y )\, &\sim\, \frac{{U} (x-v ) \wt{U} (y-u)}{n} ~
\tilde{p}_n^f (t ;v,x ;u,y ) ,
\end{split}
\end{equation}
uniformly in $x > v$, $y > u$ such that $\max\lbr |u-v|, |y-x|,|t -n\mu|\rbr \leqs \delta (n)$ 
and such that $\max\lbr x-v  ,y-u\rbr \leqs n^\epsilon$.
\medskip

\noindent
Our next task is to identify the limiting conditional (on non-intersection) marginal
distribution of the $(S_0^b ,S_1^b ,S_2^b ,\dots )$ trajectory as $n\to\infty$. The following 
notation is going to
 be convenient: Given two point $\sfw = (l ,w), \sfw^\prime  = (l ,w^\prime)\in\bbZ^2$
 with the same horizontal coordinate $l$ set
\begin{equation}
 \label{lbrconvention}
\lbr \sfw , \sfw^\prime\rbr\, \df\, (l ,w ,w^\prime )\in\bbZ^3 .
\end{equation}
Fix $m\in \bbN$. We claim:
\begin{equation}
 \label{PtildebPlus}
\begin{split}
 \lim_{n\to\infty} 
&\bbP_{v,x} \lb S_0^b = \lbr\sfv_b ,\sfx_b \rbr, S^b_1 = \lbr\sfv_1 ,\sfx_1 \rbr,
 \dots ,S^b_m = \lbr \sfv_m ,\sfx_m\rbr
~\Big| ~ S_n^{bf} = (t ,u ,y) ; \cR^{bf}_n\rb\\
&\, =\, \frac1{\wt{U} (x-v )}\bbP_{v,x}\lb 
S_0^b = (\sfv_b ,\sfx_b ), S^b_1 = \lbr\sfv_1 ,\sfx_1 \rbr, \dots ,
S^b_m = \lbr \sfv_m ,\sfx_m\rbr \rb U (x_m -v_m )\\
&\, 
\df\, 
\wt{\bbP}_{v,x}^+
\lb 
S^0_b = \lbr \sfv_b ,\sfx_b \rbr, S^b_1 = \lbr \sfv_1 ,\sfx_1 \rbr , \dots ,
S^b_m = \lbr \sfv_m ,\sfx_m\rbr \rb  , 
\end{split}
\end{equation}
as usual, 
uniformly in $x > v$, $y > u$ such that $\max\lbr |u-v|, |y-x|,|t -n\mu|\rbr \leqs \delta (n)$ 
and such that $\max\lbr x-v  ,y-u\rbr \leqs n^\epsilon$.  Indeed, formula 
\eqref{PtildebPlus} is an immediate consequence of 
\eqref{rntilda} and the second of \eqref{rntildebf}.  Notice that 
$\wt{\bbP}_{v,x}^+ $ is an instance of Doob's $h$-transform.
\medskip

\noindent
In order to develop an analogus formula for the end piece of the trajectory
and $\hat{S}_n^{b/f/bf}$ as the reversed walk, taking steps 
$\hat{\sigma}_b$, $(\hat{\sigma}_k)_{k \geq 1}$, $\hat{\sigma}_f$
(recall our definition of $\hat{\sigma}$ in \eqref{sigmahatRW}). In view of property {\bf (P3)},
satisifed by $\sigma_k$ and Remark~\ref{rem:Utilde}, $S$ and $\hat{S}$ have the same
law. On the other hand, if we set
\begin{equation}
\label{hatquantity}
 \hat{\sfu} = (0,u), \hat\sfv_b = (t -l_n, v_n), \dots,\hat\sfv_m =(t-l_{n-m} ,v_{n-m})
\end{equation}
and, accordingly, $
 \hat{\sfy} = (0,y), \hat\sfx_b = (t-l_n, x_n), \dots,\hat\sfx_m =(t- l_{n-m} ,x_{n-m})
$, then a time reversal path-transformation implies
\begin{equation}
\label{PtildefPlus}
\begin{split}
&\lim_{n\to\infty} 
\bbP_{v,x} \lb S^b_{n-m} = \lbr \sfv_{n-m} ,\sfx_{n-m} \rbr,\, S^b_{n-m+1} = \lbr \sfv_{n-m+1} ,\sfx_{n-m+1} \rbr  ,\dots, 
S^b_n = \lbr \sfv_n ,\sfx_n\rbr ~\Big| \right . \\
& \quad \quad \quad \quad  \left . ~ S_n^{bf} = (t ,u ,y) ; \cR^{bf}_n \rule[-.4cm]{0cm}{.8cm} \rb\\
&\, =\, \frac1{\wt{U} (y-u )}\bbP_{u,y}\lb 
\hat{S}^b_0 = \lbr \hat \sfv_b ,\hat \sfx_b \rbr,\, \hat{S}^b_1 = \lbr \hat \sfv_1 ,\hat \sfx_1 \rbr,\, \dots ,
\hat{S}^b_{m} = \lbr \hat\sfv_{m} ,\hat\sfx_{m} \rbr \rb U (x_{n-m} -v_{n-m })\\
&\, 
\df\, 
\wt{\bbP}_{u ,y}^+
\lb 
\hat{S}^b_0 = \lbr \hat \sfv_b ,\hat \sfx_b \rbr,\, 
\hat{S}^b_1 = \lbr \hat \sfv_1 ,\hat \sfx_1 \rbr,\, \dots ,
\hat{S}^b_m = \lbr \hat \sfv_{m} ,\hat \sfx_{m} \rbr \rb , 
\end{split}
\end{equation}
uniformly in $x > v$, $y > u$ such that $\max\lbr |u-v|, |y-x|,|t -n\mu|\rbr \leqs \delta (n)$ 
and such that $\max\lbr x-v  ,y-u\rbr \leqs n^\epsilon$. Note that under $\wt{\bbP}^+$, $S$ and $\hat{S}$ have the same distribution. Nevertheless, for the sake of clarity, we shall continue to employ them both.


\section{Repulsion and Decoupling}
\label{sec:Decoupling}
It remains to prove Lemma~\ref{lem:Decoupling}. As we have already indicated just
before the statement of the Lemma, two underlying phenomena are a repulsion 
of the trajectories of the upper and lower walks under the $\cR_n^{bf}$-constraint
 and a subsequent  asymptotic  decoupling of the event 
$\lbr \Gamma^1_\bullet\cap \gamma^{\rm up} (\Gamma_\bullet^2 ) = \emptyset \rbr$.

\noindent
With the above in mind let us proceed with 
a formal construction. First of all, repulsion will be quantified in terms of
non-intersection of certain diamond shapes.
\subsubsection{Diamond shapes} 
Given two points $\sfw$ and $\sfw^\prime$ define
 a diamond shape set
\begin{equation}
 \label{diamond}
D (\sfw , \sfw^\prime )\, =\, \lb \sfw +\cC_\delta\rb \cap\lb 
\sfw^\prime -\cC_\delta\rb .
\end{equation}
Let us say that $\lbr\Gamma^1 ,\Gamma^2\rbr\in
\cF ([\sfw ,\sfz], [\sfw^\prime , \sfz^\prime ])$ if 
\[
 \Gamma_1 = \sfw +\wt\Gamma_1,\ 
\Gamma_2 = \sfz+\wt\Gamma_2\ \ \text{and}\ \ 
\lbr\wt\Gamma^1 ,\wt\Gamma^2\rbr\in
\cF ([\sfw^\prime - \sfw ,\sfz^\prime - \sfz]) .
\]
\green{with similar definitions for $\cF_b ([\sfw ,\sfz], [\sfw^\prime , \sfz^\prime ])$
and $\cF_f ([\sfw ,\sfz], [\sfw^\prime , \sfz^\prime ])$.}
Obviously, if $\lbr\Gamma^1 ,\Gamma^2\rbr\in
\cF ([\sfw ,\sfz], [\sfw^\prime , \sfz^\prime ])$ and 
$D(\sfw ,\sfw^\prime )\cap D(\sfz , \sfz^\prime ) =\emptyset$, then
 also 
 $\Gamma^1\cap\gamma^{\rm up}(\Gamma^2 ) =\emptyset $.

\noindent
\subsubsection{The event $\cR_n^{bf}[m]$} 
Let us fix $K>0$ to be sufficiently large. 
Given an $\cR_n^{bf}$ trajectory (see the notational 
convention  \eqref{lbrconvention}),  
\[
S_0 = (0, v ,x),\, S_0^b =\lbr \sfv_b ,\sfx_b\rbr,\, S_1^b =\lbr \sfv_1 ,\sfx_1\rbr,\, \dots,\, 
S_n^b = \lbr \sfv_n ,\sfx_n \rbr ,\, S_n^{bf} = (N,u,y) ,
\]
 let us
say that $\cR_n^{bf} [m]$ occurs if, 
\begin{equation}
\label{eq:Rnplusm}
\begin{split}
&D(\sfv_k ,\sfv_{k+1})\cap D(\sfx_k ,\sfx_{k+1})\, =\, \emptyset\quad
\text{for all}\ k=m, \dots, n-m-1\\&\qquad\text{and, in addition, }\\
&l_m + (N - l_{n-m})\leq 2Km ,
\end{split}
\end{equation}
where we use the notation $\sfv_k =(l_k ,v_k)$ and, accordingly, 
$\sfx_k =(l_k ,x_k)$.
\smallskip

\noindent
Here is the key tool which enables the asymptotic analysis of
\[
\wtilde{\bbB}_p^{v,x} \lb 
\wtilde{\cA} ([\sfv ,\sfx ], [\sfu ,\sfy ])\, \big|\, 
 S_n^{bf} = (N,u,y)\, ;\,   \cR_n^{bf}\rb :
\]
\begin{prop}
\label{prop:diamonds}
There exists $\psi :\bbN\to\bbR_+$ with $\lim_{m\to\infty}\psi (m)=0$, 
such that 
\begin{equation}
\label{eq:psim}
\liminf_{N\to\infty}\bbP_{v,x}\lb \cR_n^{bf} [m]\, \big|\, 
S_n^{bf} = (N, u, y), \cR_n^{bf}\rb\,  \geq\, 1-\psi (m ) , 
\end{equation}
uniformly in $|v|, |x|, |u|, |y| \leqs \log N$ and 
$|n\mu -N|\leqs\ \sqrt{N\log N}$ .
\end{prop} 

\subsection{Repulsion}
\label{sub:RWProofs}
In this Subsection we prove Proposition~\ref{prop:diamonds}. 
Recall the we are employing the following notation for our coupled random walk: 
$S = (T,V,X)$ and $Z = X - V$.  Fix $\eta >0$ small. 
Apart from a possible violation of $T^b_m + (N-T^b_{n-m})\leq 2Km$, if 
$\cR_n^{bf} [m]$ fails to happen then either 
\[
 A_n [m]\, \df\, \lbr \exists ~k\in [m, \dots ,n-m-1]~:~ 
Z^b_k <  \min\lbr k^\eta ,(n-k)^\eta\rbr\rbr
\]
happens, or 
\[
 B_n [m]\, \df\, \lbr \exists ~k\in [m, \dots ,n-m-1]~:~\rho_k >\gamma 
\min\lbr k^\eta ,(n-k)^\eta\rbr\rbr
\]
happens, where $\gamma$ depends on the choice of the cone opening parameter 
$\delta$ in the definition \eqref{diamond} of the diamond shape $D$. In other
 words, 
\[
 \lb \cR_n^{bf} [m]\rb^{c}\, \subseteq\, \lbr T^b_m + (N-T^b_{n-m}) > 2Km\rbr
\cup A_n [m]\cup
B_n [m] .
\]
\begin{rem}
\label{rem:logparam}
Although \eqref{eq:psim} remains true for a wider range of parameters, 
all the computations will be greatly simplified if we stick to  our condition
\begin{equation}
\label{ParamRange}
|v|,|x|,|u|,|y|\leqs \log N \quad{\rm and} \quad|N-n\mu |\leqs\sqrt{N\log N} .
\end{equation}
\end{rem}
\smallskip

\subsubsection{Upper bound on $\bbP_{v,x}\lb T^b_m + (N-T^b_{n-m}) > 2Km~\big|~
\cR_n^{bf}, S_n^{bf} = (N ;u, y)\rb$}  Consider, 
\begin{equation*}
 \begin{split}
\bbP_{v, x}&\lb T^b_m >Km ,
\cR_n^{bf}, S_n^{bf} = (N ;u, y)\rb \, \\
&\leq \, \sum_{v_m < x_m}\ \sum_{s\geq Km +1}
\tilde{p}^b_m \lb s; v,x;v_m, x_m\rb \tilde{r}^f_{n-m} \lb N-s; v_m ,x_m ;u,y\rb .
\end{split}
\end{equation*}
We may ignore $|v_m | ,|x_m|,s >N^{\epsilon}$, for some $\gep < 1/4$ 
and accordingly (see \eqref{rntildebf}), 
use 
\[
\tilde{r}^f_{n-m} \lb N-s; v_m ,x_m ;u,y\rb\nnsim{1}
\frac{U(x_m -v_m )}{\wtilde{U} (x-v )}\tilde{r}_{n} \lb N; v ,x;u,y\rb .
\]
However, if $S_m = (s, v_m ,x_m)$, then (recall that we start at $(0 ,v,x)$)
\[ 
x_m -v_m\leq (x-v ) + 2\alpha s ,
\]
as it follows by the cone-confinement property {\bf (P1)} of our random walk. 
Since, in addition $r/ \wtilde{U}(r)\leqs 1$, 
we conclude that for $N$ large enough:
\begin{equation}
\label{Tmbound}
\bbP_{v,x}\lb T^b_m> Km~\big|~
\cR_n^{bf}, S_n^{bf} = (N ;u, y)\rb \, \leqs\, \bbE\left[ T^b_m ; T^b_m> Km\right] ,
\end{equation}
uniformly in $m$ and in the range of parameters described in \eqref{ParamRange}. 
The latter expression is exponentially decaying in $m$ once $K$ is fixed to be 
sufficiently large. The case of $N-T^b_{n-m} >K m$ is completely similar.\qed
\smallskip

\subsubsection{Upper bound on $\bbP_{v,x}\lb B_n [m]~\big|~
\cR_n^{bf}, S_n^{bf} = (N ;u, y)\rb$}   Write, 
\[
 B_n [m]\, =\, \bigcup_{k=m}^{n-m-1}\lbr \rho_k > \gamma\min\lbr k^\eta,  (n-k)^\eta\rbr\rbr
\]
Since time steps $\rho_k$-s have exponentially decaying tails and since by 
\eqref{rntilda} there exists $\kappa$ such that,
\[
 \tilde{r}_n (N , v,x ;u,y )\, \geqs\, \frac1{N^\kappa}
\]
uniformly in our choice of paramters in \eqref{ParamRange}, we need to consider 
only the case of $\min\lbr k^\eta,  (n-k)^\eta\rbr\leqs \log N$.  

\noindent
Again in view of exponential tails of $\rho$-variables we may 
restrict attention to $T^b_k,\, |V^b_k|,\, |X^b_k| \ll N^\epsilon$
($\gep < 1/4)$
whenever  $k^\eta \leqs \log N \nnsim{1} \log n$.
Therefore, \eqref{rntildebf} implies, 
\begin{equation*}
 \begin{split}
  \bbP_{v,x}\lb \rho_k > \gamma k^\eta ; \cR_n^{bf} ,S_n^{bf} = (N ,u, y)\rb\, &\leq\, 
\bbE_{v,x} \left[ \rho_k > \gamma k^\eta ; \tilde{r}_{n-k}^f(N- T^b_k ; V^b_k ,X^b_k; u,y)\right]\\
&\leqs\, \bbE_{v,x}\left[ \rho_k > \gamma k^\eta ;
\frac{U (X^b_k - V^b_k)}{\wtilde{U} (x-v)}\right]
\tilde{r}_n (N ; v,x ;u,y) \\
&\leqs\, 
\bbE \left[ \rho_k > \gamma k^\eta ; T^b_k\right] \tilde{r}_n (N ; v,x ;u,y)  ,
 \end{split}
\end{equation*}
where the last inequality holds for the same reason as \eqref{Tmbound}. 
Consequently, 
\[
-\log \bbP_{v,x}\lb \rho_k > \gamma k^\eta ~\big|~ \cR_n^{bf} ,
S_n^{bf} = (N ,u, y)\rb \geqs k^{\eta}
\] 
and the sum
\[ 
\sum_{k 
} \bbP_{v,x}\lb \rho_k > \gamma k^\eta ~\big|~ \cR_n^{bf}, 
S_n^{bf} = (N ,u, y)\rb
\]
converges uniformly in \eqref{ParamRange}. 
The treatment of $\bbP_{v,x}\lb \rho_k > \gamma (n-k)^\eta ~\big|~ 
\cR_n^{bf}, S_n^{bf} = (N ,u, y)\rb$ for $(n-k)^{\eta} \leqs \log N$ is similar.
\qed
\smallskip

\subsubsection{Upper bound on $\bbP_{v,x}\lb A_n [m]~\big|~
\cR_n^{bf}, S_n^{bf} = (N ;u, y)\rb$} 
As above decompose, 
\[
 A_n [m]\, =\, \bigcup_{k=m}^{n-m-1}\lbr Z^b_k < \min\lbr k^\eta,  (n-k)^\eta\rbr\rbr\, 
=\, \cup_k A_n^k .
\]
where
\[
A_n^k=
	\lbr Z^b_k < \min\lbr k^\eta,  (n-k)^\eta\rbr \rbr.
\]
Tilting, if necessary, we may assume that $\bbE T^{bf}_n = N$, and hence, 
 taking into account 
 the range of parameters in \eqref{ParamRange} and the asymptotic 
formula \eqref{rntilda},  we may assume that
\begin{equation}
\label{Anapriori}
 \tilde{r}_n (N; v,x; u,y )\, \sim\, \frac1{n^{5/2}} \wt U (x-v)\wt U (y-u ) .
\end{equation}
We shall use this as an a priori bound. Now, consider the case of $k = \min\lbr k ,n-k\rbr$:
\begin{equation*}
 \begin{split}
  \bbP_{v,x}&\lb A_n^k; \cR_n^{bf} ; S_n^{bf} = (N, u,y )\rb\\
&=\, \sum_{0 <x_k -v_k <  k^\eta} 
\sum_s \tilde{r}_k^b (s; v,x; v_k ,x_k) \tilde{r}_{n-k}^f(N-s; v_k ,x_k ;u,y ) .
 \end{split}
\end{equation*}
In view of the polynomial order of \eqref{Anapriori} it is straightforward to rule out the
possibility of (see \eqref{deltaGrowth} for properties of $\delta(\cdot )$), 
\[
\max\lbr |u- v_k |, |y-x_k|, |(N-s)- (n-k)\mu |\rbr >\delta (n/2 ). 
\]
 Hence 
for the sake of the upper bound we may assume that 
\eqref{rntildebf} uniformly  applies to all the $\tilde{r}_{n-k}$ factors above
(choose $\eta < 1/4$) 
\begin{equation}
\begin{split}
\label{rn-k}
 \tilde{r}^f_{n-k}(N-s; v_k ,x_k ;u,y )\, &\leqs\, \frac{U(x_k -v_k )\tilde U(y-u)}{n} 
\tilde{p}^f_{n-k}\lb N-s ;v_k,x_k ; u,y\rb\, \\
&\leqs\, \frac{U(k^\eta  ) \wt U (y-u)}{n} \tilde{p}_{n-k}^f 
\lb N-s ;v_k,x_k ; u,y\rb \\
&\leqs\,  \frac{U(k^\eta  ) \wt U (y-u)}{n} \tilde{p}_{n}^f 
\lb N  ;v,x ; u,y\rb
\end{split}
\end{equation}
Then, by \eqref{eqn:AprioriBdTwoVars}
\[
\begin{split}
\bbP_{v,x} 	& \lb A_n^k; \cR_n^{bf} ; S_n^{bf} = (N, u,y )\rb \\
		& \leqs \frac{U(k^\eta  ) \wt U (y-u)}{n} \tilde{p}_{n}^f \lb N  ;v,x ; u,y\rb
				\sum_{1<r<k^{\eta}} \bbP_{x,v} \lb Z^b_k=r, \cR_k^{b} \rb \\
		& \leqs	\frac{U(k^\eta  ) \wt U (y-u)}{n} \tilde{p}_{n}^f \lb N  ;v,x ; u,y\rb
				\sum_{1<r<k^{\eta}} \frac{r^2}{k^{1+1/2a^*}}
 \bbE_{v,x} \lb Z^b_0 ; Z^b_0 > 0 \rb \\
		& \leqs	\frac{U(k^\eta  ) \wt U (y-u)}{n} \tilde{p}_{n}^f \lb N  ;v,x ; u,y\rb
			(x-v) k^{3\eta - 1 - 1/2a^*}
\end{split}
\]
Hence, for $k \leq n-k$, 
\[
\bbP_{v,x} \lb A_n^k ; \cR_n^{bf} ; S_n^{bf} = (N, u,y )\rb \leqs k^{4\eta -1 - 1/2a^*}
\]

The case $n-k\leq k$ could be treated in a completely similar fashion.  It remains
to choose $\eta < 1/8$ to ensure summability of, 
\[
 \sum_{k\geq m }\frac1{k^{1 + 1/2a^* - 4\eta}}\, \sim\, \frac1{m^{1/2a^* - 4\eta }} .
\]
\qed

\subsection{Decoupling}
\noindent
\subsubsection{An a priori lower  bound} 
Define the conditional probabilities:
\begin{equation*}
\begin{split}
&p_{\sfv , \sfv^\prime}^{\sfx, \sfx^\prime} \, =\, 
\wtilde{\bbB}_p\lb \Gamma^1 \cap \gamma^{\rm up}(\Gamma^2) =\emptyset 
\, \Big| \, \lb\Gamma^1, \Gamma^2\rb \in \cF ([\sfv ,\sfx], [\sfv^\prime ,\sfx^\prime]) \rb \\
&\underline{p}_{\sfv , \sfv^\prime}^{\sfx, \sfx^\prime} \, =
\wtilde{\bbB}_p\lb \Gamma_b^1 \cap \gamma^{\rm up}(\Gamma_b^2) =\emptyset 
\, \Big| \, \lb\Gamma_b^1, \Gamma_b^2\rb \in \cF_b ([\sfv ,\sfx], [\sfv^\prime ,\sfx^\prime]) \rb \\
&\overline{p}_{\sfv , \sfu^\prime}^{\sfx, \sfx^\prime} \, =
\wtilde{\bbB}_p\lb \Gamma_f^1 \cap \gamma^{\rm up}(\Gamma_f^2) =\emptyset 
\, \Big| \, \lb\Gamma_f^1, \Gamma_f^2\rb \in \cF_f ([\sfv ,\sfx], [\sfv^\prime ,\sfx^\prime]) \rb
\end{split}
\end{equation*}

\noindent
By the finite energy property of $\otimes\Perc$ there exists $\alpha >0$, 
such that 
\begin{equation}
\label{eq:alphabound}
p_{\sfv , \sfv^\prime}^{\sfx, \sfx^\prime} \geqs \alpha^{<\sfv^\prime - \sfv, \sfe_1>}, \quad
\underline{p}_{\sfv , \sfv^\prime}^{\sfx, \sfx^\prime} \geqs \alpha^{<\sfv^\prime - \sfv, \sfe_1>}, 
\quad
\overline{p}_{\sfv , \sfv^\prime}^{\sfx, \sfx^\prime} \geqs \alpha^{<\sfv^\prime - \sfv, \sfe_1>}
\end{equation}
uniformly for all $\sfv, \sfv^\prime, \sfx, \sfx^\prime$.
On the other hand $p_{\sfv , \sfv^\prime}^{\sfx, \sfx^\prime}=1$ as soon as 
$\lbr D(\sfv ,\sfv^\prime)\cap D(\sfx ,\sfx^\prime)=\emptyset \rbr$.
In this notation the conditional $\wtilde{\bbB}_p^{v,x}$ - probability of 
$\wtilde{\cA}([\sfv , \sfx], [\sfu ,\sfy])$ 
given a trajectory 
$S^b_0 = \lbr \sfv_b ,\sfx_b\rbr,\, S^b_1 = \lbr \sfv_1 ,\sfx_1\rbr,\, \dots,\, 
S^b_n = \lbr \sfv_n ,\sfx_n\rbr,\, S^{bf}_n = \lbr \sfu ,\sfy \rbr$
from $\lbr S_n^{bf} = (N, u,y)\, ;\,   \cR_n^{bf}\rbr$
is precisely 
\[
\underline{p}_{\sfv ,\sfv_b}^{\sfx , \sfx_b}
\times {p}_{\sfv_b ,\sfv_1}^{\sfx_b , \sfx_1}
\times\prod_1^{n-1}
p_{\sfv_k , \sfv_{k+1}}^{\sfx_{k} , \sfx_{k+1}}\times 
\overline{p}_{\sfv_{n} ,\sfu}^{\sfx_{n} ,\sfy}\, \geqs \, 
\alpha^{l_b + N-l_{n} + (l_1 - l_b) + \sum_k (l_{k+1} -l_k )\1_{
\lbr D(\sfv_k ,\sfv_{k+1})\cap D(\sfx_k ,\sfx_{k+1}) \neq \emptyset \rbr} }.
\]
where we assume $\sfv_k = \lbr l_k , v_k\rbr$, $\sfx_k =\lbr l_k ,x_k \rbr$.
In view of Proposition~\ref{prop:diamonds} we infer that there
exists $\beta >0$, such that 
\begin{equation}
\label{eq:beta}
\wtilde{\bbB}_p^{v,x}\lb 
\wtilde{\cA} ([\sfv ,\sfx ], [\sfu ,\sfy ])\, \big|\, 
 S_n^{bf} = (N, u,y)\, ;\,   \cR_n^{bf}\rb\, \geq\, \beta ,
\end{equation}
uniformly in $N$, in $|v|, |x|, |u|, |y| \leqs \log N$ and in 
$|n\mu -N|\leqs\ \sqrt{N}\log N$. 

\noindent
\subsubsection{Identifying $H (\cdot )$ in  \eqref{eq:H}} 
Clearly, for every $m$,  
\begin{equation}
\label{eq:AmRm}
\wtilde{\cA}([\sfv , \sfx], [\sfu ,\sfy])\cap\cR_n^{bf}[m]\, =\, 
\wtilde{\cA}_m([\sfv , \sfx], [\sfu ,\sfy])\cap \cR_n^{bf} [m] , 
\end{equation}
where the event $\wtilde{\cA}_m =\wtilde{\cA}_m([\sfv , \sfx], [\sfu ,\sfy]) $ 
is defined exactly as event $\wtilde{\cA}$ in \eqref{eq:ATilde}, except that
the non-intersection requirement is in effect only near the boundary:
\begin{equation*}
\begin{split}
&\qquad \Gamma_b^1\cap\gamma^{\rm up}(\Gamma_b^2)=\emptyset,
\Gamma_1^1\cap\gamma^{\rm up}(\Gamma_1^2)=\emptyset,  \dots , 
\Gamma_m^1\cap\gamma^{\rm up}(\Gamma_m^2)=\emptyset\\
&\qquad\text{and, accordingly}\\
&\qquad 
\Gamma_{n-m+1}^1\cap\gamma^{\rm up}(\Gamma_{n-m+1}^2)=\emptyset,  \dots , 
\Gamma_f^1\cap\gamma^{\rm up}(\Gamma_f^2)=\emptyset .
\end{split}
\end{equation*}
Of course, $\wtilde{\cA}\subseteq\wtilde{\cA}_m$. Furthermore, 
\begin{equation}
\label{AAm}
\wtilde{\bbB}_p^{v,x`}\lb \wtilde{\cA}_m  \big|\, 
 S_n^{bf} = (N, u,y)\, ;\,   \cR_n^{bf}\rb -\, 
\wtilde{\bbB}_p^{v,x}\lb \wtilde{\cA}  \big|\, 
 S_n^{bf} = (N, u,y)\, ;\,   \cR_n^{bf}\rb\, \leq\, 2\psi (m) , 
\end{equation}
as it readily follows from \eqref{eq:psim} and \eqref{eq:AmRm}. The above bound is
uniform in $N, n, v,x,u,y$ as in the statement of 
Proposition~\ref{prop:diamonds}. In view of \eqref{eq:beta} the 
approximation is sharp (as $m\to\infty $).

\noindent
Now, conditional on a trajectory 
$\lbr \sfv_b ,\sfx_b \rbr,\, \lbr \sfv_1 ,\sfx_1 \rbr,\,\dots,\,
\lbr \sfv_n ,\sfx_n \rbr,\,\lbr \sfv_f ,\sfx_f \rbr$
the $\wtilde{\bbB}_p^{v,x}$ probability of 
$\wtilde{\cA}_m([\sfv , \sfx], [\sfu ,\sfy])$ 
is given by
\begin{equation*}
\begin{split}
&\lb \underline{p}_{\sfv ,\sfv_b}^{\sfx , \sfx_b}\times
{p}_{\sfv_b ,\sfv_1}^{\sfx_b , \sfx_1}\times
\prod_1^{m-1}
p_{\sfv_k , \sfv_{k+1}}^{\sfx_{k} , \sfx_{k+1}}\rb \, \times\,  
\lb
\prod_{n-m}^{n-1}
p_{\sfv_k , \sfv_{k+1}}^{\sfx_{k} , \sfx_{k+1}}\times
 \overline{p}_{\sfv_{n} ,\sfu}^{\sfx_{n} ,\sfy}\rb\,\\
\qquad \quad & \df\,  \underline{\frp_m} \lb \lbr \sfv,\sfx\rbr ,\lbr \sfv_b ,\sfx_b\rbr\dots \lbr \sfv_m ,\sfx_m\rbr\rb
\times \overline{\frp}_m \lb \lbr \hat{\sfu}, \hat{\sfy}\rbr ,\lbr \hat{\sfv}_b ,\hat\sfx_b\rbr ,
\dots ,\lbr \hat\sfv_m ,\hat\sfx _m\rbr\rb ,
\end{split}
\end{equation*}
where we use the same notation as in \eqref{hatquantity} (with $t=N$). 

\noindent
By the a priori bound \eqref{eq:beta}, 
\eqref{AAm}
and in view of \eqref{PtildebPlus} and \eqref{PtildefPlus}, 
 we infer that uniformly in the range of parameters \eqref{ParamRange},
\[
 \wtilde{\bbB}_p^{v,x}\lb \wtilde{\cA}_m  \big|\, 
 S_n^{bf} = (N, u,y)\, ;\,   \cR_n^{bf}\rb\lb 1 +\so \rb\, \sim\, 
 \wt{\bbE}_{v,x}^{+}\lb 
\underline{\frp_m}(S^b_0 ,\dots , S^b_m)\rb
\wt{\bbE}_{u, y}^{+}\lb 
\overline{\frp_m}(\hat S^b_0 ,\dots , \hat S^b_m)\rb 
\]
asymptotically as $n\to\infty$ and then as $m\to\infty$. 
Consequently,  we recover \eqref{eq:H}
with
\[
 H(x-v) \df \, \lim_{m\to\infty}\wt\bbE_{v,x}^{+}\lb 
\underline{\frp_m}(S^b_0 ,\dots , S^b_m)\rb
 = \lim_{m\to\infty}\wt\bbE_{v, x}^{+}\lb 
\overline{\frp_m}(\hat S^b_0 ,\dots , \hat S^b_m)\rb .
\]
\qed

\end{document}